\documentclass[11pt,reqno]{article}
\usepackage[utf8]{inputenc}
\usepackage{amsmath, mathrsfs, amssymb, amsthm, amscd, amsfonts, graphicx, color,epsfig,latexsym,eucal,layout,fancyhdr, colonequals, enumerate, hyperref, fontenc, hyperref, graphics, a4wide, verbatim, psfrag}

\newcounter{EQNR}
\renewcommand{\contentsname}{Contents\\
{\footnotesize\normalfont(A table of contents should normally not be included)}}

\newtheorem{theorem}{Theorem}
\newtheorem{corollary}[theorem]{Corollary}
\newtheorem{definition}[theorem]{Definition}

\newtheorem{lemma}[theorem]{Lemma}
\newtheorem{proposition}[theorem]{Proposition}

\newtheorem{remark}[theorem]{Remark}

\begin{document}

\title{Discrete diffusion-type equation on regular graphs and its applications
}
\author{Carlos A. Cadavid,\and Paulina Hoyos, \and Jay Jorgenson,\footnote{The third-named author acknowledges grant support
from PSC-CUNY Award 65400-00-53, which was jointly funded
by the Professional Staff Congress and The City University of New York.} \and \and \and
Lejla Smajlovi\'{c}, \and Juan D. V\'elez}
\maketitle

\begin{abstract}\noindent
We derive an explicit formula for the fundamental solution $K_{T_{q+1}}(x,x_{0};t)$ to the discrete-time diffusion equation on
the $(q+1)$-regular tree $T_{q+1}$ in terms of the discrete $I$-Bessel function.  We then use the formula to
derive an explicit expression for the fundamental solution $K_{X}(x,x_{0};t)$ to the discrete-time diffusion equation on
any $(q+1)$-regular graph $X$.
Going further, we develop three applications. The first one is to derive a general trace formula that relates the spectral
data on $X$ to its topological data.  Though we emphasize the results in the case when $X$ is finite, our method also applies
when $X$ has a countably infinite number of vertices.
As a second application,  we obtain a closed-form expression for the return time probability distribution of the uniform
random walk on any $(q+1)$-regular graph.  The expression is obtained by relating $K_{X}(x,x_{0};t)$ to the uniform random
walk on a $(q+1)$-regular graph.
We then show that if $\{X_{h}\}$ is a sequence of $(q+1)$-regular graphs whose number of vertices goes to infinity and
which satisfies a certain natural geometric condition, then the limit of the return time probability distributions from
$\{X_{h}\}$ is equal to the return time probability distribution on the tree $T_{q+1}$. As a third application, we
derive formulas which express the number of distinct closed
irreducible walks without tails on a finite graph $X$ in terms of moments of the spectrum of its adjacency matrix.
\end{abstract}

\section{Introduction}

\subsection{Diffusion-type equations}

The diffusion-type equation in timescale $T$ on the space $X$ can be viewed as the equation
\begin{equation}\label{eq: diff type eq}
\partial_t u(x;t)+ a\Delta_X u(x;t)=0,
\end{equation}
where $\partial_t$ is the partial derivative with respect to variable $t$ from the given timescale $T$ (see \cite{BP01}),
$a$ is a certain parameter denoting the strength of the diffusion (see e.g. \cite{Sl20})  and $\Delta_X $ is either the
Laplacian on the space $X$ or some generalization of it (e.g. the weighted Laplacian). The initial conditions for
\eqref{eq: diff type eq} are defined naturally and depend on the timescale.

Equation \eqref{eq: diff type eq} is rather general and there are many instances of different timescales and
different spaces $X$. The most classical situation is the diffusion equation with timescale $T=[0,\infty)$ and space
$X=\mathbb{R}$, in which case \eqref{eq: diff type eq} becomes the equation
$$
\partial_t u(x;t)= a\frac{\partial^2}{\partial x^2} u(x;t).
$$

When the space is discrete, for example when $X=\mathbb{Z}$, the diffusion-type equation \eqref{eq: diff type eq}
has been studied recently in \cite{SS14}, \cite{SS15}, \cite{Sl18}, \cite{Sl20} and \cite{Sl22}. The operator
$\Delta_X$ defined in those papers is either the centered second difference operator, which equals the
combinatorial Laplacian on the graph $\mathbb{Z}$ when viewed as a 2-regular infinite graph, or an operator
$\Delta_X$ which acts on functions $f$ defined on $\mathbb{Z}$ of the form
$$
\Delta_X f(x):= a_0f(x)-\underset{i\neq0 }{\sum_{i=-m}^{m}} a_if(x+i)
\,\,\,\,\,
\text{\rm with $x\in\mathbb{Z}$}
$$
and for some set of constants $\{a_{i}\}$. Under certain natural assumptions on the set of constants $\{a_i\}$,
such as conditions (C1)-(C4) in \cite{Sl22}, this operator can be viewed as a weighted combinatorial Laplacian on the
$2m$-regular Cayley graph with the set of vertices equal to the set of integers $\mathbb{Z}$ and the generating set being $S=\{\pm1,\ldots,\pm m\}.$

When the timescale is $T=\mathbb{N}_0$, meaning the set of non-negative integers, and $\partial_t$
is the forward difference, then the diffusion-type equation becomes a difference equation in the variable $t\in\mathbb{N}_0$, namely
\begin{equation} \label{eq: disc time general eq}
u(x;t+1)-u(x,t)+ a\Delta_X u(x;t)=0.
\end{equation}
Several instances of such an equation have been studied recently. For example,
in Theorem 4.1 of \cite{Sl18} the author considers the discrete-time diffusion equation
\begin{equation} \label{eq: Slavik u}
u(x;t+1)-u(x;t)= d(u(x+1;t) - 2u(x;t) + u(x-1;t))
\,\,\,\,\,
\text{\rm for $x\in \mathbb{Z}$}
\end{equation}
subject to the initial condition $u(x;0)=1$ if $x=0$ and $u(x;0)=0$ if $x\neq 0$. This means that in \eqref{eq: disc time general eq}
we have that $X=\mathbb{Z}$ and the operator $\Delta_{\mathbb{Z}} $ is the combinatorial Laplacian on the 2-regular graph $\mathbb{Z}$ with $a=d$.
For $d \neq 1/2$, Slav\' ik proves in \cite{Sl18} that the solution to \eqref{eq: Slavik u}, with the
 above stated initial condition, is given by
\begin{equation}\label{eq: discrete dif sol u}
u(n,t)=(1-2d)^{t}I_{|n|}^{2d/(1-2d)}(t)
\,\,\,\,\,
\text{\rm for $n\in\mathbb{Z}$ and  $t\in\mathbb{N}_0$.}
\end{equation}
Here, $I_x^c(t)$ is the discrete $I$-Bessel function introduced in \cite{Sl18} which
depends on a space variable $x\in \mathbb{Z}$, a time variable $t\in \mathbb{N}_0$, and a complex
parameter $c$; see equation \eqref{eq: I-Bessel defin} below.

The purpose of this paper is to study Equation \eqref{eq: disc time general eq} for $a=1$ when the space $X$ is an arbitrary
regular graph and $\Delta_X$ is the combinatorial Laplacian on such a graph. More precisely, we
will deduce an explicit fundamental solution to Equation \eqref{eq: disc time general eq} with
natural initial conditions and express it in terms of the graph data and the discrete $I$-Bessel functions; see Theorem \ref{thm: main} below.

We provide three applications of the main theorem. First, by studying the trace of the aforementioned fundamental solution, we deduce a trace formula for a finite regular graph; see Theorem \ref{thm: general trace f-la}. Second, by relating the combinatorial Laplacian to the stochastic Laplacian
(by which we mean the so-called random walk Laplacian), we deduce a limiting distribution result for the uniform random walk on an expanding sequence of regular graphs.
Finally, in Section \ref{sec: further counting}, by using the trace formula cited above,  we derive a formula
that expresses the number of closed irreducible walks without a tail (what are called closed geodesics with distinguished orientations in \cite{CJK14} and \cite{Mn07})
on a finite graph $X$ in terms of the moments of the spectrum of the adjacency matrix.

In the following sections we describe our main results in greater detail.

\subsection{Fundamental solution to the discrete time diffusion-type equation}

We let $X$ denote an undirected and connected graph with a finite or countably infinite set of vertices.
Let $V_{X}$ denote the set of vertices of $X$, and let $E_{X}$ denote the set of the edges of $X$.  We allow an edge
to connect a vertex to itself; such an edge is sometimes called a self-loop. We assume that $X$ is
a $(q+1)$-regular graph, where $q$ is a positive integer. This means that each vertex $x\in V_{X}$ has precisely $(q+1)$ edges that emanate
from $x$.  If $x, y \in V_{X}$, we will write $x\sim y$ if there is an edge in $E_{X}$ which connects $x$ and $y$.
The universal cover of every $(q+1)$-regular graph is a $(q+1)$-regular tree $T_{q+1}$, which is a connected
$(q+1)$-regular graph without closed walks.

Let $f:V_X \to \mathbb{C}$ be any complex-valued function on $V_{X}$.  Let $L^{2}(V_{X})$ denote the Hilbert
space of $L^{2}$ functions defined on $V_{X}$ which is  equipped with the inner product
$$
\langle f, g\rangle = \sum\limits_{x \in V_{X}}f(x)\overline{g(x)}
\,\,\,\,\,
\textrm{\rm for $f, g \in L^{2}(V_{X})$.}
$$
The Laplacian $\Delta_X$ on $X$ is a non-negative, self-adjoint operator
on $L^2(V_X)$ which is defined through its action on $L^2(V_X)$-functions, given by
$$
\Delta_X f(x)=(q+1)f(x)- \sum_{y\sim x} f(y)
\,\,\,\,\,
\textrm{\rm for $f \in L^{2}(V_{X})$.}
$$
In symbols, we can write $\Delta_X= (q+1)\mathrm{Id}-\mathcal{A}_X$, where $\mathrm{Id}$ is the identity operator
and $\mathcal{A}_X$ is the adjacency operator on the graph $X$.  The action of the adjacency operator is defined by
$$
(\mathcal{A}_Xf)(x)=\sum_{y\sim x} f(y)
\,\,\,\,\,
\text{\rm for $f \in L^{2}(V_{X})$.}
$$
If $V_{X}$ is finite, so that $X$ has a finite number of vertices say $M$, then we let
$\lambda_0=q+1 > \lambda_1\geq \ldots \geq \lambda_{M-1}\geq -(q+1)$ denote the eigenvalues of
$\mathcal{A}_X$, which are counted with multiplicities.  In this case, we let $\psi_j \in \mathbb{R}^{M}$ with $j=0,\ldots,M-1$ be the
corresponding orthonormal eigenvectors of $\mathcal{A}_X$.

Let $t\in T$ denote a time variable; we will consider either $T = \mathbb{R}_{\geq 0}$, the set of non-negative real numbers,
or $T =\mathbb{N}_0$, the set of non-negative integers.
If $T=\mathbb{R}_{\geq 0}$, then $\partial_{t,T}$ denotes the partial derivative with respect to the real variable $t$.
If $T =\mathbb{N}_0$, then $\partial_{t,T}$ denotes the forward difference operator, meaning
$$
\partial_{t,T}f(t) = f(t+1) - f(t)
\,\,\,\,\,
\text{\rm if $T =\mathbb{N}_0$.}
$$
Let $x_0\in X$ be an arbitrary but fixed point; we will call $x_{0}$ the base point of $X$.  The fundamental solution
$$
K_{X,T}: V_X\times V_X \times T \to \mathbb{R}
$$
to the diffusion-type equation
\begin{equation} \label{eq: heat eq main}
\Delta_X K_{X,T}(x_0,x;t) + \partial_{t,T} K_{X,T}(x_0,x;t) =0
\end{equation}
subject to natural initial conditions, depending on the timescale, will be called \emph{the heat kernel on the graph} $X$.
When $T=\mathbb{N}_0$, we pose the initial condition that
\begin{equation} \label{eq: heat eq initial cond}
K_{X,\mathbb{N}_0}(x_0,x;0)=\left\{
               \begin{array}{ll}
                 1, & \text{  if  } x=x_0 \\
                 0, & \text{otherwise,}
               \end{array}
             \right.
\end{equation}
while for $T=\mathbb{R}_{\geq 0}$, one replaces \eqref{eq: heat eq initial cond} by
\begin{equation} \label{eq: heat eq initial cond_cont}
\lim\limits_{t\rightarrow 0^{+}}K_{X,\mathbb{R}_{\geq 0}}(x_0,x;t)=\left\{
               \begin{array}{ll}
                 1, & \text{  if  } x=x_0 \\
                 0, & \text{otherwise.}
               \end{array}
             \right.
\end{equation}
The solution to \eqref{eq: heat eq main} subject to the initial condition
\eqref{eq: heat eq initial cond_cont} is called the \emph{continuous time heat kernel}.
When the graph $X$ is countable with bounded vertex degree, it is shown in \cite{DM06} and \cite{Do06}
that the continuous time heat kernel exists and is unique among all bounded functions.

We will refer to the solution of \eqref{eq: heat eq main} subject to the initial condition \eqref{eq: heat eq initial cond},
so $T=\mathbb{N}_0$, as the \emph{discrete time heat kernel} on the graph $X$.   When $X$ is regular of degree $(q+1)$,
the existence and uniqueness of the discrete time heat kernel follows from its connection to the uniform random walk on $X$
with initial probability distribution given by \eqref{eq: heat eq initial cond}.  This assertion is described in further detail in Lemma
\ref{lem. heat kernel in terms of random w} below.

When the time variable $t$ is continuous, there are several formulas
expressing the continuous time heat kernel on the $(q+1)$-regular tree and on a $(q+1)$-regular graph; see for example
\cite{Be97}, \cite{CY99}, \cite{CMS00}, \cite{TW03}, \cite{HNT06}, \cite{BGP09}  and \cite{CJK14}.
From the point of view taken in this article, the most interesting expression for the continuous time heat kernel is the
formula proved in \cite{CJK14}.   In that article, the authors write
the heat kernel $K_{X,\mathbb{R}_{\geq 0}}(r;t)$ in terms of the ``building blocks''
\begin{equation}\label{eq: building block cont time}
e^{-(q+1)t}q^{-r/2}\mathcal{I}_r(2\sqrt{q}t),
\end{equation}
where $r$ is the radial or distance variable between $x$ and $x_{0}$, and
as well as the set of their pre-images on $T_{q+1}$,
and $\mathcal{I}_r$ is the classical $I$-Bessel function of order $r$.

When comparing the continuous time heat kernel on $\mathbb{Z}$ from \cite{KN06} and \cite{Sl20} and the discrete time heat kernel on $\mathbb{Z}$ from \cite{Sl18},
one is led to consider the discrete $I$-Bessel functions from \cite{Sl18} to form the ``building blocks'' for the
construction of the discrete time heat kernel on the $(q+1)$-regular tree $T_{q+1}$, and
subsequently for the discrete time heat kernel on any $(q+1)$-regular graph.
Let $r\in \mathbb{N}_0$ represent the graph distance of
the point $x$ to the base point $x_0$.
We will prove below that the corresponding ``building blocks'' for the discrete time heat kernel are functions of the form
\begin{equation}\label{eq: building block disc time}
(-q)^t q^{-r/2}I_{r}^{-2/\sqrt{q}}(t)
\,\,\,\,\,
\text{\rm for $t\in\mathbb{N}_0$.}
\end{equation}
Indeed, the function $(-q)^t$ is the discrete time analogue of the
exponential function $e^{-(q+1)t}$; see Example 2.50 of \cite{BP01}. As a result,
there is a clear analogy between \eqref{eq: building block cont time} and \eqref{eq: building block disc time}.

\vskip .10in \it
For the remainder of this article, we only will consider the
discrete time case.  Hence, we will omit the notation $T$ for the timescale since we
set $T = \mathbb{N}_0$.  We will
use $K_X(x_0,x;t)$ to denote the discrete time heat kernel at points $x_{0}$ and $x$
on $X$ and $t \in \mathbb{N}_0$.
\rm
\vskip .10in

Our first main result is the following theorem; we refer the reader to section 2.1 for a detailed description of the notation and terminology which we are using.

\begin{theorem} \label{thm: main}
Let $X$ be a $(q+1)$-regular graph with a fixed base point $x_0$. For any point $x \in V_{X}$,
let $c_{m}(x)$ be the number of irreducible walks of length $m$ from $x_0$ to $x$.  Let
\begin{equation}\label{eq:b_from_c}
b_m(x)=c_m(x)-(q-1)(c_{m-2}(x)+ c_{m-4}(x)+ \ldots + c_{\ast}(x)),
\end{equation}
where
\begin{equation}\label{eq: c ast def}
c_{\ast}(x) =
\begin{cases} c_{0}(x) & \textrm{\rm if $m$ is even} \\
c_{1}(x) & \textrm{\rm if $m$ is odd}.
\end{cases}
\end{equation}
Then the discrete time heat kernel $K_X(x_0,x;t)$ on $X$ is given by
\begin{equation} \label{eq: main thm statement}
K_X(x_0,x;t)= (-q)^t\sum_{m=0}^t b_m(x) q^{-m/2}I_m^{-2/\sqrt{q}}(t)
\,\,\,\,\,
\text{\rm for $x\in V_{X}$ and $t\in \mathbb{N}_0$.}
\end{equation}
Furthermore, if we take $x=x_{0}$ and assume that $X$ is vertex transitive,
then for $x\in V_{X}$ and $t\in \mathbb{N}_0$, we have that
\begin{equation} \label{eq: main thm statement2}
K_X(x_0,x_{0};t)= (-q)^t\left(\sum_{m=0}^t N_m(x_{0}) q^{-m/2}I_m^{-2/\sqrt{q}}(t)
+(1-q)\sum\limits_{j=1}^{\lfloor t/2\rfloor}q^{-j}I_{2j}^{-2/\sqrt{q}}(t)\right)
\end{equation}
where $N_{m}(x_{0})$ is the number of closed irreducible walks from $x_{0}$ to itself of length $m$ and without tails.
\end{theorem}

The right-hand side of \eqref{eq: main thm statement} and of \eqref{eq: main thm statement2} are finite sums
because $I_m^{-2/\sqrt{q}}(t) = 0$ if $m > t$.
Hence, for any fixed $x\in V_X$, one can obtain closed-form expressions for
$b_m(x)$ as a linear combination of values of $K_X(x_0,x;t)$ at times $0\leq t \leq m$.
Also, one can obtain explicit, closed-form formulas for $c_m(x)$ in terms of $K_X(x_0,x;t)$ for $0\leq t \leq m$. This shows that geometric invariants of the graph $X$ can be determined from the discrete diffusion process on the graph in finitely many steps.
See Section \ref{sec:counting} for more details.

\subsection{Applications}

Theorem \ref{thm: main} implies the following result, which can be viewed as a type of pre-trace formula. We refer to section \ref{sec. spectral measure} for the definition of the spectral measure at the point $x$ associated to the adjacency operator $\mathcal{A}_X$ and the base point $x_{0}$.

\begin{corollary} \label{cor: pre-trace}
Let $X$ be a $(q+1)$-regular graph with a fixed base point $x_0$. Let
$\mu_x(d\lambda)$ be the spectral measure at the point $x$ associated to the adjacency operator $\mathcal{A}_X$
and the base point $x_{0}$.
Then, with the
notation from Theorem \ref{thm: main}, we have that
$$
\int\limits_{-(q+1)}^{q+1} (\lambda - q)^t \mu_x(d\lambda)=(-q)^t\sum_{m=0}^t b_m(x) q^{-m/2}I_m^{-2/\sqrt{q}}(t)
\,\,\,\,\,
\text{\rm for $x\in V_{X}$ and $t\in \mathbb{N}_0$.}
$$
When the graph $X$ is finite with $M$ vertices, then
\begin{equation}\label{eq: pre trace finite graph}
\sum_{j=0}^{M-1}(\lambda_j - q)^t \psi_j(x)\overline{\psi_j(x_0)}=
(-q)^t\sum_{m=0}^t b_m(x) q^{-m/2}I_m^{-2/\sqrt{q}}(t)
\,\,\,\,\,
\text{\rm for $x\in V_{X}$ and $t\in \mathbb{N}_0$.}
\end{equation}
\end{corollary}


When $X$ is finite, we can set $x=x_{0}$ in \eqref{eq: pre trace finite graph}
and sum over all $x \in V_{X}$ to obtain the following trace formula.

\begin{corollary} \label{cor: trace fla}
Let $X$ be a finite, undirected, $(q+1)$-regular 
graph with $M$ vertices.
Let $N_m$ denote the number of distinct closed irreducible walks without tails of length $m$ from any fixed starting point.
Then, with the notation as above, for $t \in \mathbb{N}_{0}$ we have that
\begin{equation}\label{eq: trace finite graph}
\sum_{j=0}^{M-1}\left(1-\frac{\lambda_j}{ q}\right)^t= \sum_{m=0}^t N_m q^{-m/2}I_m^{-2/\sqrt{q}}(t)
+M (1-q)\sum\limits_{j=1}^{\lfloor t/2 \rfloor}q^{-j}I_{2j}^{-2/\sqrt{q}}(t).
\end{equation}
\end{corollary}



Formula \eqref{eq: trace finite graph} can be viewed as a discrete time counterpart of the continuous time trace formula developed in \cite{Mn07}. Note that
\eqref{eq: trace finite graph} holds true for any  $(q+1)$-regular graph with $M$ vertices, while the trace formula proved in \cite{Mn07} is proved for finite connected regular graphs of degree $q+1$ without multiple edges and loops. 

In general terms, one can interpret Corollary \ref{cor: trace fla} as showing that the discrete
time heat kernel carries the same amount of spectral information about the graph as the continuous time kernel.
The advantage of the discrete time approach is that the geometric side of the expression of the heat kernel is always a
finite sum.  By comparison, we note that Theorem 1.1 of \cite{CJK14} which gives the continuous time heat kernel in terms of continuous
time $I$-Bessel functions is an infinite sum; see also \cite{Ah87} and \cite{TW03} for further results, again in continuous time.

By starting with equation \eqref{eq: trace finite graph} and using the properties of the discrete $I$-Bessel functions,
we will prove a general trace formula for any finite $(q+1)$-regular graphs.  This theorem, which we now state, is our second main result.

\begin{theorem} \label{thm: general trace f-la}
Let $X$ be a finite, undirected, $(q+1)$-regular 
graph with $M$ vertices.  For any real number $a>3+2/q$,
let $h:\mathbb{C}\to \mathbb{C}$ be a  complex valued function, holomorphic for $|z|>1/a$.  For $n \in \mathbb{N}_{0}$
and a complex number $c$, set
$$
  f_n^c(z)= \frac{1}{\sqrt{(1-z)^2 -c^2z^2} }\left(\frac{(1-z)-\sqrt{(1-z)^2 -c^2z^2}}{cz}\right)^n,
$$
where the square root is defined by taking the principal branch in the right half-plane.  Then
\begin{align*}\label{eq: general trace finite graph}
\sum_{j=0}^{M-1}h\left(\frac{q}{ q-\lambda_j} \right)&= \frac{1}{2\pi i}\sum_{m=0}^\infty N_m q^{-\tfrac{m}{2}}  \int\limits_{c(0,b)}f_m^{-2/\sqrt{q}}(z) h(z)\frac{dz}{z} \\&+ \frac{M(1-q)}{2\pi i }\sum\limits_{j=1}^{\infty}q^{-j}\int\limits_{c(0,b)}f_{2j}^{-2/\sqrt{q}}(z) h(z)\frac{dz}{z},
\end{align*}
where complex integrals above are taken along the circle centered at zero and radius $b$ for
any $b$ such that $1/a < b <q/(3q+2)$.
\end{theorem}

The heat kernel $K_X(x_0,x;t)$ is closely related to the uniform random walk on the graph $X$.
By a uniform random walk, we imagine a particle starting at an initial vertex $x_0$ at time $t=0$ and at each time step moving along one of the $(q+1)$ edges to another vertex, with all edges being selected with equal probability. Let $K_X^{\mathrm{rw}}(x_0,x;t)$ denote
the probability that the uniform random walk which starts at $x_0$ at time $t=0$ ends up at the vertex $x$ at time $t$.
It is proved in Section \ref{sec: rand walks} below that
\begin{equation}\label{eq:two_heat_kernels}
K_X^{\mathrm{rw}}(x_0,x;t)=\left(\frac{q}{q+1}\right)^t\sum\limits_{k=0}^{t}\binom{t}{k}
q^{-k}K_{X}(x_{0},x;k).
\end{equation}
Combining this result with Theorem \ref{thm: main} we obtain the following corollary.

\begin{corollary} \label{cor: rw distribution}
  Let $X$ be a $(q+1)$-regular graph. For $t\in\mathbb{Z}_{\geq 0}$, let $K_X^{\mathrm{rw}}(x_0,x;t)$
  denote the probability that the uniform random walk on the graph $X$ starts at $x_{0}$ and ends at
  $x$ at time $t$.  Then
  $$
  K_X^{\mathrm{rw}}(x_0,x;t) = \left(\frac{q}{q+1}\right)^t\sum\limits_{k=0}^{t}\binom{t}{k}
(-1)^{k} \sum\limits_{m=0}^{k}b_m(x)q^{-m/2}I_m^{-2/\sqrt{q}}(k).
  $$
\end{corollary}

When $x_0=x$, the probability distribution $K_X^{\mathrm{rw}}(x_0,x_0;t)$ equals that
of the return to $x_0$ of a random walk on $X$ after $t$ steps.  When the graph is finite,
we will compute in Section \ref{sec: limit_distribution} the limiting distribution of $K_X^{\mathrm{rw}}(x_0,x;t)$
as the number of vertices tends to infinity.  We will show that the limiting probability distribution
depends only upon $t$ and $q+1$; see \cite{TBK21}
for a related result on the first return times on finite graphs.  Further related results on $(q+1)$-regular graphs
which are Bethe lattices with coordination number $(q+1$) are given in \cite{HS82} and \cite{Gi95}, while
\cite{SRb-A05} treats the case of Erd\"os-R\'enyi random graphs.

\subsection{Organization of the paper}

The paper is organized as follows.  In Section \ref{sec:Preliminaries} we recall basic notation and recall properties of regular graphs. We prove necessary results associated to the discrete $I$-Bessel
function, including explicit formulas and asymptotic expansions in Section \ref{sec: discrete Bessel}.   Section \ref{sec: main results, proof} is devoted to proofs of our main results, except for Corollary \ref{cor: rw distribution}. The proof of  Corollary \ref{cor: rw distribution} is completed in Section \ref{sec: rand walks}. Random walks on a sequence of finite $(q+1)$-regular graphs are further studied in Section \ref{sec: limit_distribution}. We prove that the limiting distribution of return probabilities of a uniform random walk on certain sequences of finite $(q+1)$-regular graphs, as characterized by
a natural geometric condition and that the number of vertices tends to infinity, approaches the distribution of return probabilities on the tree $T_{q+1}$; see Proposition \ref{prop:limit_distribution}.  In Section \ref{sec:counting} we derive explicit and finite expressions for $N_{m}$, the number of closed irreducible walks without tails of length $m$, in terms of the spectrum of the Laplacian and special values of discrete $I$-Bessel functions.  If $X$ is finite, we show that our results complement the main results of \cite{Mn07}.  Finally, in Section \ref{sec: applications} we show how
one can solve certain other diffusion equation using the discrete $I$-Bessel function.

\subsection*{Acknowledgements}
We thank the referees of this article for their careful proofreading of a previous version of
our work.  They provided many invaluable suggestions which helped us improve
the presentation of our results.


\section{Regular graphs}\label{sec:Preliminaries}

In this section we introduce the necessary notation and properties of regular graphs.

\subsection{Basic notation}\label{sec. basic of graphs}

As above, $X$ denotes an undirected and connected graph with vertices $V_{X}$ and
edges $E_{X}$.  The set of vertices $V_{X}$ is either finite or countably infinite.  We
allow edges which connect a vertex to itself; such an edge is sometimes called a self-loop.
The degree $d_{x}$ of a vertex $x\in V_X$ is the number of vertices $y\in V_X$, counted with
multiplicity, such that
the pair $( x, y )$ is a distinct edge in $E_{X}$. The graph $X$ is \textit{vertex transitive} if its automorphism group acts transitively on the set $V_{X}$ of its vertices.

The graph $X$ is $d$-regular if $d=d_{x}$, for all $x \in V_{X}$.  Throughout this article, we only consider $d$-regular graphs.
We set $d=q+1$  and use the
parameter $q$ in our formulas.  We will choose and fix a base vertex $x_0$ in $X$.

For a non-negative integer $n$, a \textit{walk} of length $n$, as defined on page 12 in \cite{CY99},
is a sequence $y_0, y_1, \ldots, y_n$ of vertices in $V_X$ such that  $(y_j,y_{j+1})$ is an edge, with $j=0,\ldots, n-1$.
If $y_0=y_n$ the walk is a \textit{closed walk rooted at} $y_0$. A walk $y_0, y_1, \ldots, y_n$ is \textit{irreducible}
if $y_j\neq y_{j+2}$, for  $j=0,\ldots, n-2$. We say that a closed walk possesses a \textit{tail} if $y_1=y_{n-1}$.

For any two points $x,y \in V_X$ the distance between $x$ and $y$, denoted by $d_X(x,y)$, is the length of the shortest walk connecting $x$ and $y$; such a walk always exists in a connected graph. We have decided to use the word \textit{walk} instead of a
\textit{path} because some authors define paths as walks in which all adjacent vertices are distinct. In our setting we allow for loops, meaning that
our walks do not necessarily have distinct adjacent vertices.

\subsection{Spectral measure associated to the adjacency operator} \label{sec. spectral measure}

The adjacency operator $\mathcal{A}_X$ on a $(q+1)$-regular graph $X$ with a fixed base point $x_0$ is a
bounded and self-adjoint operator on $L^2(V_X)$; see  \cite[Theorem 3.1]{MW89} which describes the
results from \cite{Lo75} and \cite{Mo82}.  The spectrum $\mathrm{spec}(\mathcal{A}_X)$ of $\mathcal{A}_X$ is a subset of the interval $[-(q+1), (q+1)]$.

 According to the discussion on page 217 of \cite{MW89}, there exists a function $\mu^{\mathcal{A}_X}$ defined on the Borel sets of $\mathrm{spec}(\mathcal{A}_X)$ with values in the set of $L^2$-operators on $V_X$ such that for any continuous
 function $f$ on $\mathrm{spec}(\mathcal{A}_X)$ one has that
$$
f(\mathcal{A}_X)=\int\limits_{\mathrm{spec}(\mathcal{A}_X)} f(\lambda) \mu^{\mathcal{A}_X}(d\lambda).
$$
The function $\mu^{\mathcal{A}_X}$ is sometimes called
the resolution of the identity for $\mathcal{A}_X$.
For any point $x\in V_X$, the resolution of the identity $\mu^{\mathcal{A}_X}$ gives rise to the spectral measure $\mu_x(d\lambda)$
at $x$, associated to the adjacency operator $\mathcal{A}_X$
and the base point $x_{0}$, through the identity
$$
\mu_x(d\lambda)=\langle \mu^{\mathcal{A}_X}(d\lambda)\mathbf{e}_{x_0}, \mathbf{e}_x  \rangle
$$
where $\{\mathbf{e}_x \}_{x\in V_X}$ is the complete orthonormal system for $L^2(V_X)$.

With this notation, the spectral measure $\mu_x(d\lambda)$  at the point $x$, which is associated to the adjacency operator $\mathcal{A}_X$
and the base point $x_{0}$, has the following property.  For any continuous function $f$ defined
on $[-(q+1),(q+1)]$ one has that the evaluation of the operator $f(\mathcal{A}_X)$ at $x$ is given by
\begin{equation}\label{eq. spectral evaluation at x}
f(\mathcal{A}_X)(x)=\int\limits_{-(q+1)}^{(q+1)} f(\lambda) \mu_x(d\lambda),
\end{equation}
where we assume that the measure is identically zero outside the set $\mathrm{spec}(\mathcal{A}_X)$.

\subsection{Counting different types of walks} \label{sec: basic notation}

Let us now summarize the notation we will use which describes
certain geometric aspects of the graphs $X$ which we are studying.

\begin{definition} \label{def: graph notation} \rm Let $X$ be a $(q+1)$-regular graph.
Let $x_{0} \in X$ be a fixed vertex, and assume the notation from above.
\begin{itemize}
\item[]{1.} For any $x \in X$ and integer $m \geq 0$, let $c_{m}(x)$ be the number of irreducible walks of length $m$
from $x_0$ to $x$. In the case $x=x_{0}$, we allow closed walks with and without tails.
\item[]{2.} For any integer $m \geq 0$, let
$$
c_{m} = \sum\limits_{x\in X}c_{m}(x).
$$
\item[]{3.} For any $x \in X$, let $b_{0}(x) = c_{0}(x)=1$ and $b_{1}(x) = c_{1}(x)$.  For $m\geq 2$, let
$$
b_{m}(x) = c_{m}(x) + (1-q)\left(c_{m-2}(x) + c_{m-4}(x)+ \cdots + c_{\ast}(x) \right),
$$
where $c_{\ast}(x) $ is defined by \eqref{eq: c ast def}.
\item[]{4.} For any $x \in X$ and integer $m \geq 0$, let $N_{m}(x)$ be the number of irreducible (closed) walks of length $m$ \emph{without tails} from $x$ to $x$.
\item[]{5.} For any integer $m \geq 0$, let
$$
N_{m} = \sum\limits_{x\in X}N_{m}(x).
$$
\end{itemize}
\end{definition}

\begin{remark}\rm
The terminology employed in this article differs from that which was used in section 2.1 of \cite{CJK14}.
However, there is a correspondence which relates the terms that are used in this article and in \cite{CJK14}.
Specifically, what is called a geodesic in \cite{CJK14} is an irreducible walk, both here and in \cite{CY99}.
A closed geodesic in \cite{CJK14} means a closed irreducible walk without a tail.   We decided to
align our point of view with \cite{CY99}, which reflects the ideas and concepts of graph theory,
rather than \cite{CJK14}, which is more related to Riemannian geometry.
\end{remark}

Following the discussion in Section 2.2 of \cite{CJK14}, we obtain the following
formulas which relate the sets of values $\{N_{m}(x)\}$ and $\{c_{m}(x)\}$.

\begin{lemma}\label{lem:data_relations} Let $X$ be a $(q+1)$-regular graph.  Let $x_{0} \in X$ be a fixed vertex,
and assume the notation from above.  Additionally, assume that $X$ is vertex transitive.
\begin{itemize}
\item[]{1.} For $m \leq 2$, set $c_{m}(x_{0}) = N_{m}(x_{0})$.  Then for any $m \geq 3$, we have that
$$
N_{m}(x_{0}) = c_{m}(x_{0}) + (1-q)\left(c_{m-2}(x_{0}) + c_{m-4}(x_{0}) + \cdots + \tilde{c}_{\ast}(x_{0})\right)
$$
where
$$
\tilde{c}_{\ast}(x_{0}) =
\begin{cases} c_{2}(x_{0}) & \textrm{\rm if $m$ is even} \\
c_{1}(x_{0}) & \textrm{\rm if $m$ is odd}.
\end{cases}
$$
\item[]{2.} For any $m\geq 3$,  we have that
$$
c_{m}(x_0) =
\begin{cases} N_{m}(x_{0}) + (q-1)\sum\limits_{j=1}^{\ell-1}q^{j-1}N_{m-2j}(x_{0}) & \textrm{\rm if $m=2\ell$ is even} \\
N_{m}(x_{0}) + (q-1)\sum\limits_{j=1}^{\ell}q^{j-1}N_{m-2j}(x_{0}) & \textrm{\rm if $m=2\ell+1$ is odd}.
\end{cases}
$$
\item[]{3.} For any $m\geq 2$,  we have that
\begin{equation}\label{eq: b0 expression}
 b_{m}(x_{0})=
\begin{cases} N_{m}(x_{0}) + (1-q)  & \textrm{\rm if $m$ is even} \\
N_{m}(x_{0}) & \textrm{\rm if $m$ is odd}.
\end{cases}
\end{equation}
\end{itemize}
\end{lemma}

\begin{proof}
We will follow the proof of Proposition 2.1 of \cite{CJK14}, which itself refers \cite{Se03}.  For the sake of clarity, we will provide more detail than what is given in \cite{CJK14}.

Consider the set of irreducible walks from $x_{0}$ to itself of length $k\geq 3$ which definitely have a tail.  The number of
elements in this set is equal to $c_{k}(x_{0}) - N_{k}(x_{0})$.  Let us now provide a second approach to counting
the number of elements in this set.

Each closed, irreducible walk $\mathcal{W}_{k}$ of length $k\geq 3$ rooted at $x_{0}$ with a tail is described by a sequence of vertices which
we write as $x_{0}x_{1}x_{2}\cdots x_{k-2}x_{1}x_{0}$.
As such, the walk $x_{1}x_{2}\cdots x_{k-2}x_{1}$ can be viewed as a closed, irreducible walk $\mathcal{V}_{k-2}$ of length $k-2$ 
which is rooted at $x_{1}$ for some $x_{1}$ which is connected to $x_{0}$ by an edge.  Let us refer to the walk $\mathcal{V}_{k-2}$ as
being \it truncated \rm from $\mathcal{W}_{k}$. 


For the purposes of this lemma, it is assumed that $X$ is vertex transitive.
Hence, we can translate $\mathcal{V}_{k-2}$ to an irreducible
walk $\mathcal{W}_{k-2}$ of length $k-2$ which is rooted at $x_{0}$.  In other words, every irreducible walk $\mathcal{W}_{k}$ rooted
at $x_{0}$ can be truncated and translated to give an irreducible walk $\mathcal{W}_{k-2}$ which is also rooted at $x_{0}$.  It is immediate
that the map obtained by truncating and translating $\mathcal{W}_{k}$ to $\mathcal{W}_{k-2}$, which we denote by $\mathcal{M}_{k}(x_{0})$, is
surjective from the set of irreducible walks
of length $k$ which are rooted at $x_{0}$ and have a tail to the set of irreducible walks of length $k-2$ which are rooted at $x_{0}$.
Write $\mathcal{W}_{k-2}$ as $x_{0}y_{2}\cdots y_{k-2}x_{0}$.
Neither $y_{2}$ nor $y_{k-2}$ can equal $x_{0}$, for otherwise the walk $\mathcal{W}_{k}$ would have been reducible. 

The irreducible walks of length $k-2$ either have no tail, and there are $N_{k-2}(x_{0})$ such walks, or have a tail,
and there are $c_{k-2}(x_{0}) - N_{k-2}(x_{0})$ such walks.  It remains to determine the multiplicity of $\mathcal{M}_{k}(x_{0})$ on
each of these two disjoint sets.

If $\mathcal{W}_{k-2}$ has a tail, then it is the image of $q$ distinct walks $\mathcal{W}_{k}$.  To see this, we argue as follows.
Since $\mathcal{W}_{k-2}$ has a tail, $y_{2}=y_{k-2}$.
So, the multiplicity of the map $\mathcal{M}_{k}(x_{0})$ whose images have
tails is $q$.  Hence,
the set of all walks $\mathcal{W}_{k}$ whose image through $\mathcal{M}_{k}(x_{0})$ has a tail
has cardinality equal to $q(c_{k-2}(x_{0})-N_{k-2}(x_{0}))$. 

If $\mathcal{W}_{k-2}$ does not have a tail, then it is the image of $q-1$ distinct walks $\mathcal{W}_{k}$.  Indeed,
neither $y_{2}$ nor $y_{k-2}$ can equal $x_{0}$, 
so, the multiplicity of the map $\mathcal{M}_{k}(x_{0})$ whose images do not have
tails is $q-1$.  Hence,
the set of all walks $\mathcal{W}_{k}$ whose image through $\mathcal{M}_{k}(x_{0})$ does not have a tail
has cardinality equal to $(q-1)N_{k-2}(x_{0})$.

The set of all walks of length $k\geq 3$ with a tail can be written as a disjoint union of two sets:  One
set whose image through $\mathcal{M}_{k}(x_{0})$ has a tail, and one set whose image through
$\mathcal{M}_{k}(x_{0})$ does not have a tail.  As a result, we conclude that
\begin{equation}\label{eq:recursive_count}
c_{k}(x_{0}) - N_{k}(x_{0}) =  (q-1)N_{k-2}(x_{0}) + q(c_{k-2}(x_{0})-N_{k-2}(x_{0})),
\end{equation}
or
\begin{equation}\label{eq:recursive_countv2}
N_{k}(x_{0}) - N_{k-2}(x_{0}) = c_{k}(x_{0}) - qc_{k-2}(x_{0}).
\end{equation}
Assertion (1) follows from \eqref{eq:recursive_countv2} by induction, as does assertion (2).  Regarding (3),
the statement is a direct consequence of the definitions of $b_{m}(x_{0})$ and
$N_{m}(x_{0})$.  For odd $m$, there is no difference in the definition of
$b_{m}(x_{0})$ and $N_{m}(x_{0})$.  For even $m$, the only difference between the
two terms is that the series which defines $b_{m}(x_{0})$ ends with $c_{2}(x_{0})$ while
the series which defines $b_{m}(x_{0})$ ends with $c_{0}(x_{0})$.  Since it is immediate
that $c_{0}(x_{0})=1$, the claim follows.
\end{proof}


We now will prove a version of Lemma \ref{lem:data_relations} when $X$ is finite and not necessarily
vertex transitive.

\begin{lemma}\label{lem:data_relations2} Let $X$ be a finite $(q+1)$-regular graph with $M$ vertices.  For $m\geq 0$ let
$$
c^{(0)}_{m} = \sum\limits_{x_{0}\in X}c_{m}(x_{0}) \quad 
and \quad
b^{(0)}_{m} = \sum\limits_{x_{0}\in X}b_{m}(x_{0}).
$$
\begin{itemize}
\item[]{1.} For $m \leq 2$, set $c^{(0)}_{m} = N_{m}$.  Then for any $m \geq 3$,  we have that
$$
N_{m} = c^{(0)}_{m} + (1-q)\left(c^{(0)}_{m-2} + c^{(0)}_{m-4} + \cdots + \tilde{c}^{(0)}_{\ast}\right)
$$
where
$$
\tilde{c}_{\ast} =
\begin{cases} c^{(0)}_{2} & \textrm{\rm if $m$ is even} \\
c^{(0)}_{1} & \textrm{\rm if $m$ is odd}.
\end{cases}
$$
\item[]{2.} For any $m\geq 3$,  we have that
$$
c^{(0)}_{m} =
\begin{cases} N_{m} + (q-1)\sum\limits_{j=1}^{\ell-1}q^{j-1}N_{m-2j} & \textrm{\rm if $m=2\ell$ is even} \\
N_{m} + (q-1)\sum\limits_{j=1}^{\ell}q^{j-1}N_{m-2j} & \textrm{\rm if $m=2\ell+1$ is odd}.
\end{cases}
$$
\item[]{3.} For any $m\geq 2$,  we have that
\begin{equation}\label{eq: b0 expression2}
 b^{(0)}_{m}=
\begin{cases} N_{m} + (1-q)M  & \textrm{\rm if $m$ is even} \\
N_{m} & \textrm{\rm if $m$ is odd}.
\end{cases}
\end{equation}
\end{itemize}
\end{lemma}

\begin{proof}
Part 1. of the assertion is stated, without proof, as Proposition 2.2 in  \cite{CJK14}.  Let us provide
a proof here.  In doing so, we will use the language and notation from the proof of Lemma \ref{lem:data_relations}. 
Specifically, we will argue that the analogue of \eqref{eq:recursive_count} holds, namely that
\begin{equation}\label{eq:recursive_count2}
c^{(0)}_{k} - N_{k} =  (q-1)N_{k-2} + q(c^{(0)}_{k-2}-N_{k-2}).
\end{equation}
We view $x_{0}$ as
an arbitrary root of a closed walk, count the cardinality of certain closed walks, and then sum over all possible $x_{0}$ in order to prove \eqref{eq:recursive_count2}.  Unlike Lemma \ref{lem:data_relations}, we do not assume that $X$ is vertex transitive.

The number of closed, irreducible walks with a tail and length $k$ which are
rooted at $x_{0}$ is $c_{k}(x_{0}) - N_{k}(x_{0})$.  When summing over all $x_{0}\in X$, we have that the total number of closed,
irreducible walks with a tail and length $k$ is equal to $c^{(0)}_{k} - N_{k}$.  We will count the number of such
walks in a second manner, which will yield \eqref{eq:recursive_count2}.

As in the proof of Lemma \ref{lem:data_relations}, let $\mathcal{W}_{k}$ denote any closed, irreducible walk of length $k\geq 3$ which has a tail, described by a sequence of vertices which we write as $x_{0}x_{1}x_{2}\cdots x_{k-2}x_{1}x_{0}$.  Then $\mathcal{W}_{k}$ can be truncated to give a closed, irreducible walk $\mathcal{V}_{k-2}$ of length $k-2$ rooted at a point $x_{1}$ adjacent to $x_{0}$ and described by $x_{1}x_{2}\cdots x_{k-2}x_{1}$.  There are two cases to consider, namely when
$\mathcal{V}_{k-2}$ has a tail or when it does not have a tail. 

If $\mathcal{V}_{k-2}$ has a tail, then it is the image through truncation of $q$ different walks $\mathcal{W}_{k}$,
only one of which is rooted at $x_{0}$.  Indeed, $\mathcal{V}_{k-2}$ is the image of any walk which can be described as
$y_{0}x_{1}x_{2}\cdots x_{k-2}x_{1}y_{0}$ where $y_{0}$ is a vertex that is connected to $x_{1}$.  Since
$\mathcal{V}_{k-2}$ has a tail, $x_{2}=x_{k-2}$, so there are $q$ possibilities for $y_{0}$, namely any vertex
which is connected to $x_{1}$ by an edge and is not equal to $x_{2}$ since otherwise the extension would be reducible.
When summing over all vertices $x_{0}$, one concludes there are $q(c^{(0)}_{k-2}-N_{k-2})$ walks which are closed, irreducible
with a tail and have length $k$
whose truncations are closed, irreducible walks of length $k-2$ with a tail.

If $\mathcal{V}_{k-2}$ does not have a tail, then it is the image through truncation of $q-1$ different walks $\mathcal{W}_{k}$,
only one of which is rooted at $x_{0}$.
Indeed, $\mathcal{V}_{k-2}$ is the image of any walk which can be described as
$y_{0}x_{1}x_{2}\cdots x_{k-2}x_{1}y_{0}$ where $y_{0}$ is a vertex that is connected to $x_{1}$.  Since
$\mathcal{V}_{k-2}$ does not have a tail, $x_{2}\neq x_{k-2}$, so there are $q-1$ possibilities for $y_{0}$, namely any vertex
which is connected to $x_{1}$ by an edge and is not equal to $x_{2}$ or $x_{k-2}$ since otherwise the extension would be reducible.
When summing over all vertices $x_{0}$, one concludes there are $(q-1)N_{k-2}$ walks which are closed, irreducible
with a tail and have length $k$
whose truncations are closed, irreducible walks of length $k-2$ without a tail.

By combining the above, we arrive at equation \eqref{eq:recursive_count2}.  Then, by arguing as in the proof
of Lemma \ref{lem:data_relations}, all of the stated assertions follow.



\end{proof}
\subsection{Coverings of graphs}\label{sec. coverings of graphs}

In this section, we recall the definition of the cover of a graph from section 2 of \cite{CY99}. We will elaborate in more detail on certain aspects of sections 5  and 6 of \cite{CY99} where properties of the $(q+1)$-regular tree as a universal cover of a $(q+1)$-regular graph were discussed.

Let $\tilde{X}$ and $X$ be two simple graphs in the sense defined in \cite{CY99}, meaning that each edge has a weight equal to one. We say that $\tilde{X}$ is a \textit{covering} of $X$ if
there exists a mapping $\pi:V_{\tilde{X}} \to V_X$ which satisfies the following two properties:
\begin{itemize}
  \item [(i)] There exists a non-negative integer $m$ or possibly $m=\infty$ such that for every edge
  $(x,y)\in E_X$ the cardinality of the set $\{(\tilde{x},\tilde{y})\in E_{\tilde{X}}: \pi(\tilde{x})=x,\, \pi(\tilde{y})=y\}$ equals $m$.
  \item [(ii)] For $\tilde{x},\, \tilde{y} \in V_{\tilde{X}}$ with $\pi(\tilde{x})=\pi(\tilde{y})$ and $y$ adjacent to $\pi(\tilde{x})$ in $X$, the
  sets $N_{\tilde{X}}(\tilde{x}) \cap \pi^{-1}(y)$ and  $N_{\tilde{X}}(\tilde{y}) \cap \pi^{-1}(x)$ have the same cardinality, where
  $N_{\tilde{X}}(\tilde{x})$ denotes the
  set of all vertices in $V_{\tilde{X}}$ adjacent to $\tilde{x}$ in $\tilde{X}$.
\end{itemize}

When $\tilde{X}=T_{q+1}$ is the $(q+1)$-regular tree and $X$ is a $(q+1)$-regular graph, then there exists a natural mapping $\pi: V_{T_{q+1}} \to V_X $ such that for each $\tilde{x} \in T_{q+1}$ the neighbors of $\tilde{x}$ are mapped to  neighbors of $\pi(\tilde{x})$ in $X$ in one-to one way and the root $\tilde{x}_0$ of the tree $T_{q+1}$ is mapped to the fixed base vertex $x_0$ of $X$.

For a fixed integer $r\geq 0$, and any $x\in V_X$, it is elementary to observe that the number of elements $\tilde{x}$ in $\pi^{-1}(x)$ such that the distance $d_{T_{q+1}}(\tilde{x}_0,\tilde{x})$ between the root $\tilde{x}_0$ of the tree and $\tilde{x}$ equals the number of irreducible walks in $X$ from $x_0$ to $x$ of length $r$.  

\section{Discrete $I$-Bessel functions}\label{sec: discrete Bessel}

Let $c\in \mathbb{C}$ be an arbitrary complex number, and let $t,n \in\mathbb{N}_0$
be non-negative integers.  By expanding ideas from \cite{BC18}, the author of \cite{Sl18}
defines the discrete modified $I$-Bessel function $I_n^c(t)$ as follows.
Let
$$
(-t)_0=1
\,\,\,\,\,
\text{\rm and}
\,\,\,\,\,
(-t)_n=(-t)(-t+1)\cdot \ldots \cdot (-t+n-1)
\,\,\,\,\,
\text{\rm for $n>0$}
$$
be the Pochhammer symbol.  Recall that the (Gauss) hypergeometric series is defined by
$$
F(\alpha,\beta;\gamma;z):= \sum_{k=0}^{\infty} \frac{(\alpha)_k(\beta)_k}{k! (\gamma)_k} z^k;
$$
see section 9.1 of \cite{GR07}.  Then, the discrete $I$-Bessel function $I_n^c(t)$ is defined by
\begin{equation} \label{eq: I-Bessel defin}
I_n^c(t) := \frac{(-c/2)^n (-t)_n}{n!}F\left(\frac{n-t}{2},\frac{n-t}{2}+\frac{1}{2}; n+1;c^2 \right).
\end{equation}

For any $n$ and $m \in \mathbb{N}_{0}$, we have that
\begin{equation}\label{eq:Poch_vanishing}
(-m)_n=0
\,\,\,\,\,\text{\rm for}
\,\,\,\,\,
n>m.
\end{equation}
Hence,  $I_n^c(t)=0$ for $n>t$.

\subsection{Basic properties of the discrete $I$-Bessel function}

In the lemma below we summarize some properties of the discrete $I$-Bessel function $I_n^c(t)$ which are proved in \cite{Sl18}
and will be needed in this article.  As above, $\partial_{t}$ denotes the forward difference operator, meaning
that $\partial_{t}f(t) = f(t+1)-f(t)$.

\begin{lemma}  \label{lem: properties of discrete I-Bes}
With the notation as above, the discrete $I$-Bessel function $I_n^c(t)$ has the following properties.
\begin{itemize}
  \item[(i)] $I_n^c(t)=0$ for $n>t$.
  \item[(ii)] $I_0^c(0)=1$.
  \item[(iii)] $\partial_t I_1^c(0)=c/2$ and $\partial_t I_0^c(t)=cI_1^c(t)$ for $t\geq 1$.
  \item[(iv)] $\partial_t I_x^c(t)=\frac{c}{2}\left(I_{x-1}^c(t) + I_{x+1}^c(t)\right)$ for any $x \geq 1$ and $t\geq 0$.
\end{itemize}
\end{lemma}

If $n\leq t$, one of the numbers $\frac{n-t}{2}$ and $\frac{n-t}{2}+\frac{1}{2}$ must be a non-positive integer.
As a result, by \eqref{eq:Poch_vanishing}, the hypergeometric series in formula \eqref{eq: I-Bessel defin} yields a
polynomial in the variable $c^2$ whenever $n \leq t$.  Therefore, $I_n^c(t)$ is a polynomial in the variable $c$.  The
following proposition gives an explicit evaluation of $I_n^c(t)$ in this instance.

\begin{proposition} \label{prop: I as polynomial}
Let $t,n \in\mathbb{N}_0$ such that $n\leq t$.  Set $\ell =\lfloor (t-n)/2\rfloor$.
Then for any $c\in \mathbb{C}$, we have that
\begin{equation}\label{eq: I bessel as polynom}
  I^{c}_{n}(t)=\sum\limits_{j =0}^{ \ell }\frac{t!}{j
!(t-2j -n)!(n+j )!}\left(\frac{c}{2}\right) ^{2j +n}.
\end{equation}
\end{proposition}

\begin{proof}
For any $d\in\mathbb{R}\setminus\{ 1/2\}$, consider the discrete time diffusion equation \eqref{eq: Slavik u}
subject to the initial condition $u(x;0)=1$ if $x=0$ and $u(x;0)=0$ for $x \neq 0$.
From Theorem 4.1 of \cite{Sl18}, we have that the solution
for this equation, with the restriction that $d\neq 1/2$, is given by
\eqref{eq: discrete dif sol u}.

On the other hand, from Example 3.3 of \cite{SS14}, we also have that
$$
u(n,t)=(1-2d)^t\sum_{j=0}^t \binom{t}{j,t-2j-n,j+n}\left(\frac{d}{1-2d}\right)^{2j+n};
$$
see in addition Remark 4.2 of \cite{Sl18}.  (Note: We follow
the common convention that multinomial coefficients with negative terms are equal to zero.)
Therefore, when comparing the two expressions for $u(n,t)$, we deduce that
$$
  I^{c}_{n}(t)=\sum_{j=0}^t \binom{t}{j,t-2j-n,j+n}\left(\frac{c}{2}\right)^{2j+n}
$$
for $n\in\mathbb{N}_0$ and any real number $c \neq -1$.
Since $t-2j-n <0$ for $j> \lfloor(t-n)/2 \rfloor$, the assertion is proved for
all $c\in\mathbb{R}\setminus\{-1\}$.  From this, note that one can view
$I^{c}_{n}(t)$ as a function of the complex variable $c$.  As such, the proposition follows
for all $c\in\mathbb{C}$ by analytic continuation.
\end{proof}

From Proposition \ref{prop: I as polynomial}, we have that $I_n^c(t)$ is a polynomial in the variable $c$ of degree
$2\lfloor \frac{t-n}{2} \rfloor +n$.  In other words, $I_n^c(t)$ is a polynomial in variable $c$ of degree $t$ when $t-n$ is even
and of degree $t-1$ when $t-n$ is odd.

Going further, it is necessary to determine the asymptotic behavior of the ``building block'' $(-q)^t q^{-r/2}I_{r}^{-2/\sqrt{q}}(t)$
as $q\to \infty$.   To this end, we have the following result.

\begin{corollary}
For positive integers $q,r,t$ such that $r\leq t$ we have that
$$
(-q)^t q^{-r/2}I_{r}^{-2/\sqrt{q}}(t) = (-1)^{t-r}\binom{t}{r}q^{t-r}+ O(q^{t-r-1})
\,\,\,\,\,\,\,\,
\text{\rm as $q \rightarrow \infty$.}
$$
\end{corollary}

\begin{proof}
  For a fixed $r,t$ such that $r\leq t$, equation \eqref{eq: I bessel as polynom} with $c=-2/\sqrt{q}$ and $\ell=\lfloor(t-r)/2 \rfloor$
  gives that
  \begin{align*}
    (-q)^t q^{-r/2}I_{r}^{-2/\sqrt{q}}(t) &= (-q)^t q^{-r/2} \sum\limits_{j =0}^{ \ell }\frac{t!}{j!(t-2j -r)!(r+j )!}\left(\frac{-1}{\sqrt{q}}\right) ^{2j +r} \\&= (-1)^{t+r}q^{t-r} \sum\limits_{j =0}^{ \ell }\frac{t!}{j!(t-2j -r)!(r+j )!} q^{-j}
    \\&= (-1)^{t-r}\binom{t}{r}q^{t-r} + O(q^{t-r-1})
  \end{align*}
as $q\to\infty$, which completes the proof.
\end{proof}

\begin{remark} \rm
  Let $t,n \in\mathbb{N}_0$ such that $n\leq t$ and put $\ell =\lfloor (t-n)/2\rfloor$.
  Recall that the Jacobi polynomial $P_{n}^{\alpha,\beta}(x)$ can be expressed in terms of
  the hypergeometric function, and also we have that $P_{n}^{\alpha,\beta}(x)$ has a precise
  expression whose coefficients are given in terms of combinatorial coefficients; see 8.962.1
  and 8.960 of \cite{GR07}, respectively.  By combining those expressions, we can derive another
  formula for the discrete $I$-Bessel function.  Namely, for any $c\in \mathbb{C}$, we have that
\begin{equation}\label{eq: I bessel as polynom 1}
I_n^c(t)=\left(\frac{c}{2}\right)^n \frac{\binom{t}{n}}{\binom{n+\ell}{n}}\sum_{m=0}^{\ell}
\binom{\ell - (-1)^{t-n}/2}{m}\binom{\ell + n}{\ell - m} c^{2m}.
\end{equation}
We find the representation \eqref{eq: I bessel as polynom} to be more desirable since it is
reminiscent of the series representation of the classical $I$-Bessel function, which is
$$
\mathcal{I}_n(z)=\sum_{k=0}^{\infty}\frac{1}{k!(n+k)!}\left(\frac{z}{2}\right)^{2k+n}.
$$
\end{remark}

\subsection{Asymptotic behavior of discrete $I$-Bessel functions}\label{sec: Asympt I-Bessel}

In this section, we will determine the asymptotic behavior of the discrete $I$-Bessel function $I_n^c(t)$
for a fixed real parameter $c$ and as $t\to \infty$. For our analysis, we will use the representation
\eqref{eq: I-Bessel defin} of the function $I_n^c(t)$, recalling that the asymptotics of the hypergeometric
function for the large values of parameters is well studied; see for example \cite{Te03}.  Specifically, we will use
the transformations of the hypergeometric function given in formulas 9.122, 9.131, and 9.132 from \cite{GR07}.
From these formulas, the asymptotic behavior of the discrete $I$-Bessel function is obtained by considering the
different cases studied in \cite{Wa18} and \cite{Jo01} where the authors proved certain
asymptotic and uniform asymptotic results.

\begin{proposition}
For any real, nonzero parameter $c$ and a fixed $n\in \mathbb{N}$, we have
that
\begin{equation}\label{eq: I asympt}
I_n^c(t) \sim \frac{(\mathrm{sgn}(c))^n}{\sqrt{2\pi |c| t}}(1+|c|)^{t+1/2}
\,\,\,\,\,
\text{\rm as $t \rightarrow \infty$,}
\end{equation}
where $\mathrm{sgn}(c)$ denotes the sign of $c$.
\end{proposition}

\begin{proof}
First consider the case when $\vert c \vert = 1$.  From
  formula 9.122 of \cite{GR07} we have that
  $$
  F\left(\frac{n-t}{2},\frac{n-t}{2}+\frac{1}{2}; n+1;1 \right)=\frac{\Gamma(n+1)\Gamma\left(t +\frac{1}{2}\right)}{\Gamma\left(\frac{n+t+1}{2}\right)\Gamma\left(\frac{n+t+1}{2} +\frac{1}{2}\right)}.
  $$
  The classical duplication formula for the gamma function and definition \eqref{eq: I-Bessel defin} yields that
  $$
  I_n^{\pm1}(t)= (\pm 1)^n\cdot 2^{-n}\frac{t(t-1)\cdot \ldots \cdot (t-(n-1))}{n!}\cdot \frac{n!\Gamma\left(t +\frac{1}{2}\right) 2^{t+n}}{\sqrt{\pi}\Gamma\left(t+n+1\right)}.
  $$
  Using the functional equation for the gamma function we then get that
   $$
  I_n^{\pm1}(t)= (\pm 1)^n\cdot \frac{t(t-1)\cdot \ldots \cdot (t-(n-1))}{(t+n)(t+n-1)\cdot \ldots \cdot (t+1)}\cdot \frac{\Gamma\left(t +\frac{1}{2}\right) 2^{t}}{\sqrt{\pi}\Gamma\left(t+1\right)}\sim (\pm 1)^n \frac{2^t}{\sqrt{\pi t}}
  $$
  as $t \to \infty$. This proves \eqref{eq: I asympt} when $c\in\{-1,1\}.$

Next, assume that $|c|<1$.  To begin, we apply formula 9.131 of \cite{GR07} which gives that
$$
F\left(\frac{n-t}{2},\frac{n-t}{2}+\frac{1}{2}; n+1;c^2 \right)=(1-c^2)^{t/2-n/2}F\left(\frac{t+n+1}{2},\frac{n-t}{2}; n+1;\frac{c^2}{c^2-1} \right).
$$
Since $|c|<1$, we can write $\frac{c^2}{c^2-1}=\frac{1}{2}(1-z)$ for $z=\frac{1+c^2}{1-c^2}>1$.  Hence,
$z+\sqrt{z^2-1}= \frac{1+|c|}{1-|c|}=e^{\zeta}$, for $\zeta>0$, in the notation of page 289 of \cite{Wa18}.
By applying the asymptotic formula from \cite{Wa18}, we get that
$$
F\left(\frac{t+n+1}{2},\frac{n-t}{2}; n+1;\frac{c^2}{c^2-1} \right)\sim \frac{2^n}{\sqrt{\pi t}} \frac{n!\Gamma\left(\frac{t-n}{2} +1\right)}{\Gamma\left(\frac{t+n}{2} +1\right)} \left( \frac{2|c|}{1+|c|} \right)^{-n-\tfrac{1}{2}}
\left( \frac{1+|c|}{1-|c|}\right)^{\tfrac{t}{2}-\tfrac{n}{2}}
$$
as $t \rightarrow \infty$.  Now, let us
multiply the above asymptotic formula by $\frac{(-c/2)^n (-t)_n}{n!}(1-c^2)^{t/2-n/2}$ and use that
$$
\frac{(-1)^n (-t)_n}{2^n}\frac{\Gamma\left(\frac{t-n}{2} +1\right)}{\Gamma\left(\frac{t+n}{2} +1\right)} \sim 1
$$
as $t\to\infty$.  In doing so, we deduce the claimed asymptotic formula for $|c|<1$,
thus proving \eqref{eq: I asympt} for $|c|<1$.

Finally, assume that $|c|>1$.  In this case, we proceed analogously as above.
In the application of the asymptotic formulas from \cite{Wa18}, the lead term in the absolute value is
$\frac{1+|c|}{1-|c|}=e^{-\zeta}$ where $\zeta=x+i\pi$ with $x<0$. Therefore, we have that
\begin{align*}
F&\left(\frac{t+n+1}{2},\frac{n-t}{2}; n+1;\frac{c^2}{c^2-1} \right) \\ &\sim \frac{2^n}{\sqrt{\pi t}} \frac{n!\Gamma\left(\frac{t-n}{2} +1\right)}{\Gamma\left(\frac{t+n}{2} +1\right)} \left( \frac{2|c|}{|c|-1} \right)^{-n-\tfrac{1}{2}}(-1)^{\tfrac{n+1}{2}}\left( \frac{1+|c|}{1-|c|}\right)^{\tfrac{t}{2}+\tfrac{n+1}{2}}
\,\,\,\,\,
\text{\rm as $t \rightarrow \infty$.}
\end{align*}
One then multiplies this asymptotic formula by $\frac{(-c/2)^n (-t)_n}{n!}(1-c^2)^{t/2-n/2}$ and, by arguing as
above,  we obtain \eqref{eq: I asympt}.  With all this, the proof of the assertion is complete.
\end{proof}

The above proposition is very useful in finding the limiting behavior of solutions to the heat equation when the time variable tends to infinity. A simple consequence of formulas \eqref{eq: discrete dif sol u} and \eqref{eq: I asympt} is the following corollary:
\begin{corollary}
  For $d>0$, consider the solution \eqref{eq: discrete dif sol u} to the discrete diffusion equation \eqref{eq: Slavik u} subject to the initial condition $u(x;0)=1$ if $x=0$ and $u(x;0)=0$ if $x \neq 0$.  Then, one has that
  $$
  u(n,t)\sim
\frac{1}{2\sqrt{\pi d t}}
\,\,\,\,\,
\text{\rm for $d\in (0,1/2)$ and as $t\rightarrow \infty$}
$$
and
$$u(n,t)\sim
\frac{(-1)^{n+t}(4d-1)^{t+1/2}}{2\sqrt{\pi d t}}
\,\,\,\,\,
\text{\rm for $d\in (1/2,\infty)$ and as $t\rightarrow \infty$.}
$$
\end{corollary}

\subsection{A generating function of discrete $I$-Bessel functions}\label{sec: gen f-on I-Bessel}

For now, let $c \in\mathbb{C}\setminus \{0\}$  and  $n\in\mathbb{N}_0$. Since $|I_n^c(t)|\leq I_n^{|c|}(t)$ for all integers $n,t\geq 0$, the asymptotic formula \eqref{eq: I asympt} implies that the radius of convergence of the power series
\begin{equation}\label{eq: gen funct of I}
f^c_n(z):= \sum_{t=0}^{\infty} I_n^c(t)z^t
\end{equation}
is at least $(1+|c|)^{-1}$.  Therefore, the series \eqref{eq: gen funct of I} is a holomorphic function
of $z\in \mathbb{C}$ in the disc $|z|<(1+|c|)^{-1}$.
We will now derive a closed expression for the generating function \eqref{eq: gen funct of I}
when $|z|<(1+|c|)^{-1}$.

\begin{proposition} \label{prop: gen funct}
  For any $c\in\mathbb{C}\setminus \{0\}$, $n\in\mathbb{N}_0$, and $|z|<(1+|c|)^{-1}$,
  the power series \eqref{eq: gen funct of I} is equal to
  $$
  f_n^c(z)= \frac{1}{\sqrt{(1-z)^2 -c^2z^2} }\left(\frac{(1-z)-\sqrt{(1-z)^2 -c^2z^2}}{cz}\right)^n,
  $$
where we are using the principal value of the square root, which takes values in the right half-plane.
\end{proposition}

\begin{proof}
 We start by observing that $$(1-z)f_n^c(z)=I_n^c(0)+\sum\limits_{t=0}^{\infty} \left(I_n^c(t+1)-I_n^c(t)\right)z^{t+1}= I_n^c(0)+z\sum\limits_{t=0}^{\infty}\partial_tI_n^c(t)z^t.$$
 Using all four statements of Lemma \ref{lem: properties of discrete I-Bes}, one can show that the power series \eqref{eq: gen funct of I}
  satisfies the difference equation
\begin{equation}\label{eq: recurence for fc}
  (1-z)f_n^c(z)=\left\{
                       \begin{array}{ll}
                         1+czf_{n+1}^c(z), & \text{  if  } n=0; \\
                         \frac{cz}{2}\left(f_{n-1}^c(z) + f_{n+1}^c(z)\right), & \text{  if  } n\geq 1.
                       \end{array}
                     \right.
\end{equation}
To solve \eqref{eq: recurence for fc}, let us set the Ansatz that $f_n^c(z) =a(z)(b(z))^n$, where we have omitted, for convenience,
$c$ from the  notation on the right-hand side.
When $n=0$,  \eqref{eq: recurence for fc} yields that $a(z)=(1-z-czb(z))^{-1}$.  In the case $n\geq 1$,
equation \eqref{eq: recurence for fc} implies
that
$$
b(z)^2 - \frac{2(1-z)}{cz}b(z)+1=0
$$
for $|z|<(1+|c|)^{-1}$ and $z\neq 0$.  There are two solutions to this equation, namely
$$
b_{1}(z)= \frac{1-z}{cz}+ \sqrt{\frac{(1-z)^2 -c^2z^2}{c^2z^2}}
\,\,\,\,\,\text{\rm and}\,\,\,\,\,
b_{2}(z)= \frac{1-z}{cz}- \sqrt{\frac{(1-z)^2 -c^2z^2}{c^2z^2}}.
$$
In order to choose the solution which yields the generating function, recall that $f_n^c(0)=I_n^c(0)$.
Furthermore, $I_0^c(0)=1$ and $I_n^c(0)=0$ for $n \geq 1$.
Therefore, the function $b(z)$ must be such that $\lim_{z\to 0}b(z)=0$, and this is true when
$$
b(z) = b_{2}(z)= \frac{(1-z)- \sqrt{(1-z)^2 -c^2z^2}}{cz}.
$$
Therefore, we have that $a(z)=( \sqrt{(1-z)^2 -c^2z^2})^{-1}$, and the proof of the assertion is complete.
\end{proof}

\section{Proofs of main results}\label{sec: main results, proof}

In this section, we prove our main results. We start by proving an expression for the heat kernel on a $(q+1)$-regular tree
in terms of the discrete $I$-Bessel function.  As stated above, the formula is analogous to the result from \cite{CJK14}
which gives the continuous time heat kernel on the $(q+1)$-regular tree in terms of the continuous $I$-Bessel function.

\subsection{The discrete time heat kernel on a $(q+1)$-regular tree}

Let us fix a base point $x_0\in X$. Let $T_{q+1}$ be a $(q+1)$-regular tree, which is the universal cover of the graph $X$ with base point $x_0$, as described in section \ref{sec. coverings of graphs}.  Choose a lift of $x_{0}\in X$ to $T_{q+1}$ to give a base point of $T_{q+1}$
which in a slight abuse of notation we also denote by $x_{0}$ (this point was denoted by $\tilde{x}_0$ in section \ref{sec. coverings of graphs}).
Let $r\in \mathbb{N}_0$ be the radial coordinate on $T_{q+1}$ with the chosen base point $x_{0}$.
Recall that the radial coordinate of a point $x$ is the graph distance $d_{T_{q+1}}(x_0,x)$ from $x$ to the base point $x_0$, see section \ref{sec. basic of graphs}.

As before, $t\in \mathbb{N}_0$ and $\partial_t$ is the forward difference operator with respect to the variable $t$.
Let $K_{q+1}(x_{0},x;t)$ be the discrete time heat kernel on $T_{q+1}$.  There is an automorphism $g$ of $T_{q+1}$ to
itself, which fixes $x_{0}$ and sends any point $x'$ with radial coordinate $r$ to any other point $x$ with radial
coordinate $r$. It is immediate, using the characterizing properties \eqref{eq: heat eq main} and
\eqref{eq: heat eq initial cond} (with $X=T_{q+1}$ in \eqref{eq: heat eq main}) that $K_{q+1}(gx_{0},gx;t)=K_{q+1}(x_{0},x;t)$.
 Therefore, the heat kernel $K_{q+1}(x_{0},x;t)$ can be viewed as a function of $r$ and $t$, meaning
$$
K_{q+1}(x_{0},x;t) = K_{q+1}(r;t)
\,\,\,\,\,
\text{\rm where $r = r(x)= d_{T_{q+1}}(x_0,x)$.}
$$

In order to derive the characterizing equations for the heat kernel on $T_{q+1}$ in radial coordinates, one should notice that the root $x_0$ of $T_{q+1}$ is adjacent to $q+1$ vertices, all of which are at radial distance $1$ from $x_0$. On the other hand, if $x\neq x_0$ is at radial distance $r\geq 1$ from $x_0$, then it is adjacent to $q$ vertices which are at radial distance $r+1$ from $x_0$ and to one vertex that is at radial distance $r-1$ from $x_0$. Therefore, the action of $\Delta_{T_{q+1}}= (q+1)\mathrm{Id} - \mathcal{A}_{T_{q+1}}$ on $K_{q+1}(x_{0},x;t) = K_{q+1}(r;t)$ is, in the radial coordinate $r$, given by
$$
\Delta_{T_{q+1}} K_{q+1}(r;t) =\left\{
                       \begin{array}{ll}
                         (q+1)K_{q+1}(0;t) - (q+1)K_{q+1}(1;t), & \text{  if  } r=0 \\
                         (q+1)K_{q+1}(r;t)-qK_{q+1}(r+1;t)-K_{q+1}(r-1;t), & \text{  if  } r\geq 1.
                       \end{array}
                     \right.
$$
Therefore, the characterizing properties \eqref{eq: heat eq main} and \eqref{eq: heat eq initial cond}, when
given in radial coordinates, reduce to the equation
\begin{equation}\label{eq: heat eq on a tree}
  \partial_t K_{q+1}(r;t)= \left\{
                       \begin{array}{ll}
                         -(q+1)K_{q+1}(0;t) + (q+1)K_{q+1}(1;t), & \text{  if  } r=0 \\
                         -(q+1)K_{q+1}(r;t)+qK_{q+1}(r+1;t)+K_{q+1}(r-1;t), & \text{  if  } r\geq 1,
                       \end{array}
                     \right.
\end{equation}
with initial condition
\begin{equation}\label{eq: initial cond tree}
  K_{q+1}(r;0)= \left\{
            \begin{array}{ll}
              1, & \text{  if  } r=0 \\
              0, & \text{  if  } r\geq 1.
            \end{array}
          \right.
\end{equation}

\begin{proposition} \label{Prop: HK on tree}
With the notation as above, the solution $K_{q+1}(r;t)$
to the initial value problem \eqref{eq: heat eq on a tree} and \eqref{eq: initial cond tree} is given by
\begin{equation}\label{eq: heat on the tree expresion}
  K_{q+1}(r;t)=(-q)^t\left[ q^{-r/2}I_r^{-2/\sqrt{q}}(t) - (q-1) \sum_{j=1}^{\infty} q^{-(r+2j)/2}I_{r+2j}^{-2/\sqrt{q}}(t)\right].
\end{equation}
\end{proposition}

\begin{proof}
Parts (i) and (ii) of Lemma \ref{lem: properties of discrete I-Bes} imply that the right-hand side of
\eqref{eq: heat on the tree expresion} satisfies the initial condition \eqref{eq: initial cond tree} when $t=0$.
Hence, it is left to prove that the right-hand side of \eqref{eq: heat on the tree expresion} satisfies
\eqref{eq: heat eq on a tree}.  To begin, let us define
$$
F(r;t):= q^{-r/2}I_r^{-2/\sqrt{q}}(t) - (q-1) \sum_{j=1}^{\infty} q^{-(r+2j)/2}I_{r+2j}^{-2/\sqrt{q}}(t).
$$
If the above stated assertion is true, then $F(r;t) = (-q)^{-t} K_{q+1}(r;t)$.
It is elementary to verify that the right-hand side of \eqref{eq: heat on the tree expresion}
satisfies \eqref{eq: heat eq on a tree} if and only if $F(r;t)$ satisfies the equation
\begin{equation}\label{eq: heat for F}
   \partial_t F(r;t)= \left\{
                       \begin{array}{ll}
                          - \frac{(q+1)}{q}F(1;t), & \text{  if  } r=0 \\
                         -F(r+1;t)-\frac{1}{q}F(r-1;t), & \text{  if  } r\geq 1.
                       \end{array}
                     \right.
\end{equation}
We now shall prove that \eqref{eq: heat for F} is true.
When $r=0$ we have that
$$
\partial_t F(0;t)= \partial_tI_0^{-2/\sqrt{q}}(t)-(q-1) \sum_{j=1}^{\infty} q^{-j}\partial_tI_{2j}^{-2/\sqrt{q}}(t).
$$
Using parts (iii) and (iv) of Lemma \ref{lem: properties of discrete I-Bes} we get that
\begin{align*}
\partial_t F(0;t) &= -\frac{2}{\sqrt{q}}I_1^{-2/\sqrt{q}}(t)-(q-1)\sum_{j=1}^{\infty} q^{-j}\frac{-2}{2\sqrt{q}}\left(I_{2j+1}^{-2/\sqrt{q}}(t)+I_{2j-1}^{-2/\sqrt{q}}(t)\right)\\&=
I_1^{-2/\sqrt{q}}(t)\left( -\frac{2}{\sqrt{q}} + \frac{q-1}{q\sqrt{q}}\right) + (q-1)\sum_{j=1}^{\infty} \frac{q^{-(j+1)}(q+1)}{\sqrt{q}}I_{2j+1}^{-2/\sqrt{q}}(t)
\\&=-\frac{q+1}{q}\left(q^{-1/2}I_1^{-2/\sqrt{q}}(t) - (q-1) \sum_{j=1}^{\infty} q^{-(2j+1)/2}I_{2j+1}^{-2/\sqrt{q}}(t)\right)
\\&=-\frac{q+1}{q} F(1;t),
\end{align*}
which proves \eqref{eq: heat for F} in the case when $r=0$.

Assume now that $r\geq 1$.  From part (iv) of Lemma \ref{lem: properties of discrete I-Bes} we get that
\begin{align*}
\partial_t F(r;t) &= -\frac{2q^{-r/2}}{2\sqrt{q}}\left(I_{r+1}^{-2/\sqrt{q}}(t)+ I_{r-1}^{-2/\sqrt{q}}(t)\right)
\\&\hskip .25in -(q-1)\sum_{j=1}^{\infty} q^{-(r+2j)/2}\frac{-2}{2\sqrt{q}}\left(I_{r+2j+1}^{-2/\sqrt{q}}(t)+I_{r+2j-1}^{-2/\sqrt{q}}(t)\right)\\&=
-q^{-(r+1)/2}I_{r+1}^{-2/\sqrt{q}}(t)+ (q-1)\sum_{j=1}^{\infty} q^{-(r+1+2j)/2}I_{r+2j+1}^{-2/\sqrt{q}}(t)\\&\hskip .25in -\frac{1}{q}q^{-(r-1)/2}I_{r-1}^{-2/\sqrt{q}}(t) + \frac{q-1}{q}\sum_{j=1}^{\infty} q^{-(r-1+2j)/2}I_{r+2j-1}^{-2/\sqrt{q}}(t)\\&=
-\left(F(r+1;t)+\frac{1}{q}F(r-1;t)\right).
\end{align*}
This proves \eqref{eq: heat for F} when $r\geq 1$, and completes the proof of the proposition.
\end{proof}

\begin{remark}\rm
Note that for fixed $t\in \mathbb{N}_0$, the series on the right-hand side of formula \eqref{eq: heat on the tree expresion}
is actually a finite sum.  Indeed, $I_{r+2j}^{-2/\sqrt{q}}(t)=0$ for all $j$ such that $r+2j>t$, so the terms
in the series in \eqref{eq: heat on the tree expresion} are non-zero only if
$j$ is such that $r+2j \leq t$.
\end{remark}

\subsection{Proof of Theorem \ref{thm: main}}

Let $X$ be any $(q+1)$-regular graph.  Let $T_{q+1}$ be the universal cover of $X$ with covering map
$\pi$, as described in section \ref{sec. coverings of graphs}. Then we can write the heat kernel $K_X(x_0,x;t)$ on $X$ by
\begin{equation}\label{eq:heat_periodize}
K_X(x_0,x;t) = \sum\limits_{\tilde{x}\in \pi^{-1}(x)}K_{q+1}(x_{0},\tilde{x};t).
\end{equation}

The natural covering map $\pi$ has the property that it maps the immediate neighbours of $\tilde{x}\in T_{q+1}$ to neighbours of $\pi(\tilde{x})$,
and this mapping is injective.  Therefore, points $\tilde{x}_0=x_0$ and $\tilde{x}\in\pi^{-1}(x)$ have
distance $r$ in $T_{q+1}$ if and only if $x_0$ and $x$ can be connected through an irreducible walk of
length $r$. Let $c_r(x)$ be the number of irreducible walks in $X$ of length $r$ from $x_0$ to $x$. Then, $K_{q+1}(x_{0},\tilde{x};t)=K_{q+1}(r;t)$  for exactly $c_r(x)$ elements $\tilde{x}\in \pi^{-1}(x)$, hence \eqref{eq:heat_periodize} can be written as
$$
K_X(x;t)=\sum_{r\geq 0} c_r(x)K_{q+1}(r;t).
$$
Since the base point $x_{0}$ fixed, we will suppress that portion of the notation
and simply write $K_X(x;t) = K_X(x_0,x;t)$.

Note that the series on the right-hand side of the above display is actually a finite sum for any fixed
$t\in \mathbb{N}_0$ since $K_{q+1}(r;t)=0$ when $r>t$.  As a result, there is no question
regarding the convergence of the series.

When using the expression \eqref{eq: heat on the tree expresion} for $K_{q+1}(r;t)$, we immediately obtain that
$$
K_X(x;t)= (-q)^t \sum_{r\geq 0} c_r(x)\left(q^{-r/2}I_r^{-2/\sqrt{q}}(t) - (q-1) \sum_{j=1}^{\lfloor\frac{t-r}{2} \rfloor}
q^{-(r+2j)/2}I_{r+2j}^{-2/\sqrt{q}}(t) \right).
$$
For a fixed positive integer $t$, and fixed $r\leq t$, let us choose $r_1\leq t$. Then, the term $I_{r_1}^{-2/\sqrt{q}}(t)$
on the right hand side of the above equation arises when $r_1=r$ and for all values of $j$ in the set $\{1,\ldots, \lfloor\frac{t-r}{2}
\rfloor\}$ for which $r+2j=r_1$. Therefore, the factor multiplying $I_{r_1}^{-2/\sqrt{q}}(t)$ equals
$$
(-q)^t\left( c_{r_1}(x)q^{-r_1/2} - (q-1)q^{-r_1/2} \left(c_{r_1-2}(x)+ c_{r_1-4}(x)\ldots c_{\ast}(x)\right) \right)
$$
where $c_{\ast}(x)$ equals $c_0(x)$ if $r_1$ is even and $c_1(x)$ if $r_1$ is odd. This proves  \eqref{eq: main thm statement}.

Equation \eqref{eq: main thm statement2} follows by combining \eqref{eq: main thm statement} and \eqref{eq: b0 expression}.

\subsection{Proofs of corollaries \ref{cor: pre-trace} and \ref{cor: trace fla}}

Let us write \eqref{eq: heat eq main} as
\begin{align}\nonumber
K_{X}(x_0,x;t+1)&=(\mathcal{A}_X-q\mathrm{Id})K_{X}(x_0,x;t) \\&= (\mathcal{A}_X-q\mathrm{Id})^t K_{X}(x_0,x;0).\label{eq. heat k in terms of adjacency}
\end{align}
When we apply formula \eqref{eq. spectral evaluation at x} with $f(\mathcal{A}_X)= (\mathcal{A}_X-q\mathrm{Id})^t$, we deduce that
$$
K_{X}(x_0,x;t)=\int\limits_{-(q+1)}^{q+1} (\lambda - q)^t \mu_x(d\lambda).
$$
With this, the first part of Corollary \ref{cor: pre-trace} is proved.

When $X$ is a finite graph of degree $q+1$ with $M$ vertices, the adjacency operator can be identified with the adjacency
matrix $A_X$ which possesses real eigenvalues $\lambda_0=q+1 > \lambda_1\geq \ldots \geq \lambda_{M-1}\geq -(q+1)$, counted
according to their multiplicities.  The associated orthonormal eigenvectors are denoted by $\psi_j \in \mathbb{R}^{M}$ for $j=0,\ldots,M-1$.
In this case, the spectral expansion of the heat kernel reads as
$$
K_{X}(x_0,x;t)=\sum_{j=0}^{M-1}(\lambda_j - q)^t \psi_j(x)\overline{\psi_j(x_0)},
$$
which proves the second part of Corollary \ref{cor: pre-trace}.

It remains to prove Corollary \ref{cor: trace fla}. We first take $x=x_0$ in \eqref{eq: pre trace finite graph}.
After multiplying \eqref{eq: pre trace finite graph} by $(-q)^{-t}$, we then get that
\begin{equation} \label{eq: heat kernel finite M}
\sum_{j=0}^{M-1}\left(1-\frac{\lambda_j}{ q}\right)^t\psi_j(x_{0})\overline{\psi_j(x_{0})}=
 \sum_{m=0}^t b_{m}(x_{0}) q^{-m/2}I_m^{-2/\sqrt{q}}(t) .
\end{equation}
Recall that the eigenvectors $\psi_j$ are normalized to have $L^2$-norm equal to $1$.
Now, by taking the sum over all vertices $x_{0} \in V_X$ and using Lemma \ref{lem:data_relations2}, part 3,
the result follows.

\subsection{Proof of Theorem \ref{thm: general trace f-la}}

Let $h(z)$ be a function which is holomorphic for $\vert z \vert > 1/a$ for some positive real number $a$.
Let us assume that $a>3+2/q > (2q+1)/q$.
As such, $h(z)$  can be written as a Taylor series centered at $z_0=\infty$, namely
$$
h(z)=\sum_{t=0}^{\infty} g(t) z^{-t}.
$$
In this notation,  $h(z)$ is the one-sided $\mathcal{Z}$-transform of function $g:\mathbb{N}_0\to \mathbb{C}$.
Equivalently, we can say that $\{g(t)\}$ is the set of Taylor series coefficients of $h$.
Moreover, the convergence of the series which defines $h$, together
with assumption that $a>3+2/q > (2q+1)/q$ implies the existence of some $\epsilon > 0$ with
$0<\epsilon<\frac{1}{2}(a-(2q+1)/q)$ such that
\begin{equation}\label{eq:g_bound}
g(t)((2q+1)/q +\epsilon)^t \to 0
\,\,\,\,\,
\text{\rm as $t \rightarrow \infty$.}
\end{equation}

Since $aq>2q+1$ and  $-(q+1)\leq \lambda_j\leq q+1$ for every eigenvalue $\lambda_{j}$, we have that
$\left|\frac{q}{q-\lambda_j}\right|>1/a$. Therefore,
\begin{equation}\label{eq:lhs_trace_formula}
\sum_{j=0}^{M-1}h\left(\frac{q}{q-\lambda_j}\right) = \sum_{t=0}^{\infty} g(t)\sum_{j=0}^{M-1} \left(1-\frac{\lambda_j}{q}\right)^t.
\end{equation}
The above expression can be obtained from the left-hand side of \eqref{eq: trace finite graph}
after multiplying \eqref{eq: trace finite graph} by $g(t)$ and then summing over all $t\geq 0$.
Let us now apply the same operations to the right-hand side of \eqref{eq: trace finite graph}
and compute the outcome.

As stated, $(2q+1)/q \geq 1+2/\sqrt{q}$ for $q\geq 1$.  In \eqref{eq: I asympt} we derived the
asymptotic behavior of the discrete $I$-Bessel function as $t\to \infty$, which we now will use
with $c=-2/\sqrt{q}$.  By combining \eqref{eq: I asympt} with \eqref{eq:g_bound},
we conclude that the series $\sum_{t=0}^{\infty} g(t)  I_m^{-2/\sqrt{q}}(t)$ converges absolutely.

Moreover, from Proposition \ref{prop: gen funct} and the power series representation \eqref{eq: gen funct of I}, for all $z$ such that $1/a<|z|<(1+2/\sqrt{q})^{-1}$, we have that
\begin{equation}\label{eq:fh_product}
f_n^{-2/\sqrt{q}}(z)h(z)= \sum_{t_1=0}^{\infty}I_n^{-2/\sqrt{q}}(t_1)z^{t_1} \sum_{t_2=0}^{\infty}g(t_2)z^{-t_2}.
\end{equation}
Note that the annulus $1/a<|z|<(1+2/\sqrt{q})^{-1}$ is not empty because we assumed that $a>(2q+1)/q$
and $(2q+1)/q \geq 1+2/\sqrt{q}$.
The residue theorem implies that
the constant term in \eqref{eq:fh_product} can be expressed in terms of the integral of $f_n^{-2/\sqrt{q}}(z)h(z) z^{-1}$
along any circle inside the annulus $1/a<|z|<(1+2/\sqrt{q})^{-1}$.  Specifically, for any $b$ such
that $1/a<b<q/(3q+2)$ we have that
\begin{equation}\label{eq: product expression}
\sum_{t=0}^{\infty}I_n^{-2/\sqrt{q}}(t)g(t)= \frac{1}{2\pi i} \int\limits_{c(0,b)}f_n^{-2/\sqrt{q}}(z) h(z)\frac{dz}{z}.
\end{equation}
 In order to complete the proof, it suffices to show that the series
\begin{equation}\label{eq: final series}
\sum_{m=0}^\infty N_m q^{-\tfrac{m}{2}} \frac{1}{2\pi} \int\limits_{c(0,b)}\left|f_m^{-2/\sqrt{q}}(z) h(z)\frac{dz}{z}\right|
\end{equation}
is convergent. In order to do so, we estimate function $f_m^{-2/\sqrt{q}}(z)$ for $z$ with $|z|=b<q/(3q+2)$. Recall that
$$
  f_m^c(z)= \frac{1}{\sqrt{(1-z)^2 -c^2z^2} }\left(\frac{cz}{(1-z)+\sqrt{(1-z)^2 -c^2z^2}}\right)^m.
$$
Trivially, we have that
$$
\vert (1-z)+\sqrt{(1-z)^2 -c^2z^2}\vert \geq \text{\rm Re} ( (1-z)+\sqrt{(1-z)^2 -c^2z^2})\geq \text{\rm Re}(1-z)\geq (1-b).
$$
Moreover, for $q=4$, the function $(1-z)^2 -4q^{-1}z^2$ vanishes at $z=1/2$, while for $q\neq 4$ it vanishes at $z=\frac{q\pm 2\sqrt{q}}{q-4}$. It is straightforward to see that $\left|\frac{q\pm 2\sqrt{q}}{q-4}\right| >\frac{1}{2}>b$ for all $q\geq 1$ provided $q\neq 4$.
Hence, for any $q \geq 1$, we have that $(1-z)^2 -4q^{-1}z^2$ is non-vanishing on the circle $|z|=b$. Therefore,
$$
\vert f_m^{-2/\sqrt{q}}(z)\vert \leq M(b,q) \left(\frac{2bq^{-1/2}}{1-b}\right)^m
$$
where $M(b,q)$ is the maximal value of the continuous function $|\sqrt{(1-z)^2 - 4q^{-1}z^2}|^{-1}$ on the circle $|z|=b$.
The function $h(z)$ is holomorphic,
 hence bounded by a constant $C(h,b)$ along the circle $|z|=b$. The number $N_m$ of closed irreducible walks without a
 tail of length $m$ is bounded by the number of total walks of length $m$, so then $N_m\leq (q+1)^m$.  Thus,
$$
 N_m q^{-\tfrac{m}{2}} \frac{1}{2\pi} \int\limits_{c(0,b)}\left|f_m^{-2/\sqrt{q}}(z) h(z)\frac{dz}{z}\right|\leq M(b,c) C(h,b) \left(\frac{2b(q+1)}{q(1-b)}\right)^m.
$$
The value of $b$ is chosen so that $b<q/(3q+2)$, hence $\frac{2b(q+1)}{q(1-b)}<1$.  Therefore,
the series \eqref{eq: final series} is bounded by a convergent geometric series.

With all this, we can summarize the proof of Theorem \ref{thm: general trace f-la} as follows.
First, multiply equation \eqref{eq: trace finite graph} by $g(t)$ and take the sum over all $t\geq 0$.
Then, interchange the summation over $m$ and $t$ and arrive at the equation
$$
\sum_{j=0}^{M-1}h\left(\frac{q}{q-\lambda_j}\right)= \sum_{m=0}^{\infty}N_m q^{-m/2} \sum_{t=0}^{\infty} g(t)
I_m^{-2/\sqrt{q}}(t) +M(q-1)\sum_{j=0}^{\infty}q^{-2j} \sum_{t=0}^{\infty} g(t)  I_{2j}^{-2/\sqrt{q}}(t).
$$
Finally, by inserting the identity \eqref{eq: product expression}, one obtains Theorem \ref{thm: general trace f-la}.
In effect, the detailed analysis above justifies the operations described which yield Theorem \ref{thm: general trace f-la}
from \eqref{eq: trace finite graph}.

\begin{remark}\rm
Corollary \ref{cor: rw distribution} follows directly from Theorem \ref{thm: main}
and equation \eqref{eq:two_heat_kernels}.  The proof of \eqref{eq:two_heat_kernels}
is given in section \ref{sec: rand walks} below.  Once the proof of \eqref{eq:two_heat_kernels}
is established, all statements in the introduction are proven.
\end{remark}

\section{Uniform random walk on a $(q+1)$-regular graph} \label{sec: rand walks}

Let $X$ denote a $(q+1)$-regular graph, which is either finite or countably infinite, and which has a base point $x_0$.
A uniform random walk on $X$ starting at $x_0$ is the discrete time random walk at which the particle located at any
vertex of $X$ at time $t=n\geq 0$ moves to any of its $(q+1)$ neighboring vertices.  The movement to
any vertex occurs with probability $1/(q+1)$.

This process can be described as a solution to the diffusion equation associated to the uniform random walk Laplacian
$\Delta_X^{\mathrm{rw}}$ on $X$ which is given by
$$
\Delta_X^{\mathrm{rw}}=\mathrm{Id}-\frac{1}{q+1}\mathcal{A}_{X}
= \frac{1}{q+1}\Delta_X.
$$

Specifically, the probability that a particle which starts walking at the base point $x_0$ and after $t$ steps is located at the point $x$ equals
$ K_{X}^{\mathrm{rw}}(x_0,x;t)$ where $ K_{X}^{\mathrm{rw}}(x_0,x;t)$ is the discrete time random walk heat kernel.  In other words,
$K_{X}^{\mathrm{rw}}(x_0,x;t)$  is the
solution to the discrete time diffusion equation
\begin{equation} \label{eq: heat RW eq main}
\Delta_X^{\mathrm{rw}} K_{X}^{\mathrm{rw}}(x_0,x;t) + \partial_t K_{X}^{\mathrm{rw}}(x_0,x;t) =0
\end{equation}
subject to the initial condition
\begin{equation} \label{eq: heat RW eq initial cond}
 K_{X}^{\mathrm{rw}}(x_0,x;0)=\left\{
               \begin{array}{ll}
                 1, & \text{  if  } x=x_0 \\
                 0, & \text{otherwise.}
               \end{array}
             \right.
\end{equation}

Note that \eqref{eq: heat RW eq main} can be written as
\begin{equation}\label{eq. heat rw in terms of adjacency}
 K_{X}^{\mathrm{rw}}(x_0,x;t)=\frac{1}{q+1}\mathcal{A}_{X}K_{X}^{\mathrm{rw}}(x_0,x;t)=
 \frac{1}{(q+1)^{t}}\mathcal{A}_{X}^{t}K_{X}^{\mathrm{rw}}(x_{0},x;0),
\end{equation}
which succinctly describes the dynamics of this process.

The heat kernels $K_{X}(x_0,x;t)$ and $K_{X}^{\mathrm{rw}}(x_0,x;t)$ have the same initial condition, meaning that
\begin{equation}\label{eq:equal_intitial_conditions}
K_{X}^{\mathrm{rw}}(x_{0},x,0)=K_{X}(x_{0},x,0)
\,\,\,\,\,
\text{\rm for all $x$.}
\end{equation}
With this, we can solve for one heat kernel in terms of the other. The result is given in the following lemma.
\begin{lemma} \label{lem. heat kernel in terms of random w}
 Let $X$ be a $(q+1)$-regular graph, which is either finite or countably infinite, with a base point $x_0$. For any $x\in X$ the discrete time heat kernel $K_{X}(x_0,x;t)$ and the discrete time random walk heat kernel $K_{X}^{\mathrm{rw}}(x_0,x;t)$ are related through the identities
\begin{equation} \label{eq: relation stanrard and rw hk}
K_{X}(x_0,x;t)=(-q)^t\sum_{k=0}^t(-1)^k \binom{t}{k}\left(1+\frac{1}{q}\right)^k   K_{X}^{\mathrm{rw}}(x_0,x;k).
\end{equation}
and \eqref{eq:two_heat_kernels}, which is
$$
K_X^{\mathrm{rw}}(x_0,x;t)=\left(\frac{q}{q+1}\right)^t\sum\limits_{k=0}^{t}\binom{t}{k}
q^{-k}K_{X}(x_{0},x;k).
$$
\end{lemma}
\begin{proof}
To prove \eqref{eq: relation stanrard and rw hk} we combine \eqref{eq. heat k in terms of adjacency},
\eqref{eq:equal_intitial_conditions} and \eqref{eq. heat rw in terms of adjacency}.  Specifically, we have that
\begin{align*}
K_{X}(x_0,x;t) & = (\mathcal{A}_{X}-q\mathrm{Id})^{t}K_{X}(x_{0},x;0) \\&= (\mathcal{A}_{X}-q\mathrm{Id})^{t} K_{X}^{\mathrm{rw}}(x_{0},x;0)
\\&= \sum\limits_{k=0}^{t}\binom{t}{k}(-q)^{t-k}(\mathcal{A}_{X})^{k}K_{X}^{\mathrm{rw}}(x_{0},x;0)
\\&= (-q)^{t}\sum\limits_{k=0}^{t}(-1)^{k}\binom{t}{k}\left(\frac{q+1}{q}\right)^{k}
\frac{1}{(q+1)^{k}}(\mathcal{A}_{X})^{k}K_{X}^{\mathrm{rw}}(x_{0},x;0)
\\&= (-q)^{t}\sum\limits_{k=0}^{t}(-1)^{k}\binom{t}{k}\left(1+\frac{1}{q}\right)^{k}
K_{X}^{\mathrm{rw}}(x_{0},x;k),
\end{align*}
as claimed.  On the other hand, for $\mathcal{B}_{X}= \mathcal{A}_{X}-q\mathrm{Id}$, one has that
\begin{eqnarray*}
K_{X}^{\mathrm{rw}}(x_{0},x;t) &=&\frac{1}{(q+1)^{t}}(\mathcal{B}_{X}+q\mathrm{Id)%
}^{t}K_{X}^{\mathrm{rw}}(x_{0},x;0)\\ \nonumber
&=&\frac{1}{(q+1)^{t}}\sum\limits_{k=0}^{t}\binom{t}{k}q^{t-k}\mathcal{B}%
_{X}^{k}K_{X}(x_{0},x;0) \\
&=&\frac{1}{(q+1)^{t}}\sum\limits_{k=0}^{t}\binom{t}{k}%
q^{t-k}K_{X}(x_{0},x;k),
\end{eqnarray*}
which proves \eqref{eq:two_heat_kernels}.
\end{proof}


By combining \eqref{eq:two_heat_kernels} with \eqref{eq: main thm statement}, we have proved Corollary \ref{cor: rw distribution}.

As a special case of the above computations, one can take $X$ to be the
$(q+1)$-regular tree.  By combining \eqref{eq: heat on the tree expresion} with  \eqref{eq:two_heat_kernels}, we immediately obtain
an explicit expression of the random walk heat kernel $K_{q+1}^{\mathrm{rw}}(r;t)$ associated to a $(q+1)$-regular tree in terms of the discrete
$I$-Bessel function.  Specifically, for $t \geq r$, one has that
\begin{equation}\label{eq: rw heat kernel f-la}
K_{q+1}^{\mathrm{rw}}(r;t)=\frac{q^t}{(q+1)^t}\sum_{k=r}^t \binom{t}{k}(-1)^k \left(q^{-r/2}I_r^{-2/\sqrt{q}}(k) - (q-1) \sum_{j=1}^{\lfloor \frac{k-r}{2}\rfloor} q^{-(r+2j)/2}I_{r+2j}^{-2/\sqrt{q}}(k)\right),
\end{equation}
where, as before, $r$ denotes the radial variable on the $(q+1)$-regular tree relative to a fixed base point $x_{0}$.

\begin{remark}\rm
The random walk heat kernel on homogeneous regular trees was studied in \cite{Ur97}, Section 2, where certain properties of
$K_{q+1}^{\mathrm{rw}}(r;t)$ were derived. However, those results were based on recurrence formula satisfied by this
kernel and no closed formula for its evaluation was deduced. A more general setting of semi-regular infinite graphs
was studied in \cite{Ur03} where the upper and lower bounds for the discrete time random walk heat kernel were deduced.  Again,
no closed formula for this heat kernel was obtained.
\end{remark}

\begin{remark}\rm
By taking $r=0$ in equation \eqref{eq: rw heat kernel f-la} we can deduce the return probability after $t$ steps of a random walk on a
$(q+1)$-regular tree, meaning the probability that a uniform random walk on the tree starting at the root comes back to the root
after $t$ steps. This probability is given by
\begin{equation}\label{eq: rw return prob tree}
K_{q+1}^{\mathrm{rw}}(0;t)=\frac{q^t}{(q+1)^t}\sum_{k=0}^t \binom{t}{k}(-1)^k \left(I_0^{-2/\sqrt{q}}(k) - (q-1) \sum_{j=1}^{\lfloor \frac{k}{2}\rfloor} q^{-j}I_{2j}^{-2/\sqrt{q}}(k)\right).
\end{equation}
\end{remark}

\section{A limiting distribution result} \label{sec: limit_distribution}

The quantity $K_{X}^{\mathrm{rw}}(x_0,x_0;t)$ can be viewed as the
return probability of the random walk on $X$ to the starting point after $t$ steps.
Using the above formulas, we can determine the limiting behavior of the return
probability for certain sequences of finite graphs as the number of vertices of $M$ tends to infinity.

\begin{proposition}\label{prop:limit_distribution}
  Let $\{X_{h}\}$ be a sequence of finite $(q+1)$-regular vertex transitive graphs with base points $\{x_{0,h}\}$.  Assume that
  $X_{h}$ has $M_{h}$ vertices, and that $M_{h}$ tends to infinity as $h$ tends to
  infinity.  Assume further that the length of the shortest closed irreducible walk without a tail on $X_{h}$ rooted at $x_{0,h}$ tends to infinity as $M_{h}$ goes to infinity.
  Then for fixed $t \geq 0$ we have, for sufficiently large $h$, the equality
  \begin{equation} \label{eq: limit of rw }
  K_{X_{h}}^{\mathrm{rw}}(x_{0,h},x_{0,h};t)= \frac{q^t}{(q+1)^t}
  \sum_{k=0}^t\binom{t}{k}(-1)^{k}\left(I_0^{-2/\sqrt{q}}(k) -(q-1)\sum\limits_{j=1}^{\lfloor \frac{t}{2}\rfloor}q^{-j}
  I_{2j}^{-2/\sqrt{q}}(k)\right).
  \end{equation}
\end{proposition}

\begin{remark}\rm In a slight abuse of notation, we will write the statement in Proposition
\ref{prop:limit_distribution} as saying that for fixed $t$, we have that
$$
  K_{X}^{\mathrm{rw}}(x_{0,h},x_{0,h};t)= K_{q+1}^{\mathrm{rw}}(0;t).
$$
for $h$ sufficiently large.
\end{remark}

\begin{proof}
  Let $X$ be any fixed $(q+1)$-regular graph with $M$ vertices. Then, from \eqref{eq: heat kernel finite M} and \eqref{eq: b0 expression},
  we have that the heat kernel $K_X(x_{0,h},x_{0,h};t)$ can be written as
  $$
  K_X(x_{0,h},x_{0,h};t) = (-q)^t\left[ \sum_{m=0}^t N_m (x_{0,h}) q^{-m/2}I_m^{-2/\sqrt{q}}(t) -
  (q-1)\sum\limits_{j=1}^{\lfloor \frac{t}{2}\rfloor}q^{-j}I_{2j}^{-2/\sqrt{q}}(t)\right],
  $$
  where, as stated above, $N_m(x_{0,h})$ is the number of closed irreducible walks without a tail of length $m$ beginning at $x_{0,h}$.
  For the sequence $\{X_{h}\}$ under consideration, and for fixed $t$, it is assumed that for sufficiently large $h$
  we have that $N_{m}(x_{0,h}) = 0$ for all $0<m\leq t$.
  Recall the convention that $N_0(x_{0,h})=1$.  Therefore, for fixed $t$ and sufficiently large $h$, we have that
  \begin{equation}\label{eq:limit_return_time}
   K_X(x_{0,h},x_{0,h};t) = (-q)^{t}\left(I_0^{-2/\sqrt{q}}(t) -
   (q-1)\sum\limits_{j=1}^{\lfloor \frac{t}{2}\rfloor}q^{-j}I_{2j}^{-2/\sqrt{q}}(t)\right).
  \end{equation}

By combining this with formula \eqref{eq:two_heat_kernels}, the proof is complete.
\end{proof}

\begin{remark}\rm
According to \cite{TBK21}, the limit as $M \to\infty$ of the \it first return probability distribution \rm of a random
walk on a $(q+1)$-regular randomly chosen graph with $M$ vertices equals the first return probability distribution
on the $(q+1)$-regular tree.  The above computations show that one has a similar result when considering the (total) return
probability distribution for certain sequences of graphs since the return probability distribution \eqref{eq: rw return prob tree} on the tree $T_{q+1}$ coincides
with the distribution obtained by combining \eqref{eq:limit_return_time} with \eqref{eq:two_heat_kernels}.
\end{remark}

\section{Relating spectral data to length spectrum data}\label{sec:counting}

As with trace formulas in general, one can use Theorem \ref{thm: general trace f-la} to express
spectral data of the Laplacian, or adjacency operator, to topological data, namely the length spectrum.
In the case $X$ is finite, such expressions are stated in Corollary 2 of \cite{Mn07}, under the additional assumption that $X$ is connected, without multiple edges and loops.  In this
section, we will use Corollary \ref{cor: trace fla} to show how to express the length spectrum, meaning
the number $N_m$ of distinct closed irreducible walks without tails of length $m$,
in terms of moments of the spectrum of the adjacency matrix $A_X$. Going further, we will
discuss in Section \ref{sec: further counting} below how to use Corollary
\ref{cor: pre-trace} to prove similar results in the case $X$ is an infinite $(q+1)$-regular graph.

\subsection{Explicit evaluations of $N_{j}$ for $j \leq 3$}
Let us explicitly compute $N_{j}$ for all $j \leq 3$ for any finite $(q+1)$-regular graph with $M$ vertices
using formula \eqref{eq: trace finite graph}.  In each case, the evaluations use the precise formula for
the discrete $I$-Bessel function, as stated in \eqref{eq: I bessel as polynom}.

Using that $I_{0}^{-2/\sqrt{q}}(0)=1$ and $I_{2}^{-2/\sqrt{q}}(0)=0$, we get that equation \eqref{eq: trace finite graph}
in the case $t=0$ yields the formula $N_{0}=M$, as is set by convention.

Let us now consider $t=1$.  From Lemma \ref{lem: properties of discrete I-Bes},
we can readily compute that $I_{0}^{-2/\sqrt{q}}(1) = 1$, $I_{1}^{-2/\sqrt{q}}(1) = -2/(2\sqrt{q})$ and $I_{2}^{-2/\sqrt{q}}(1) = 0$.
Thus, equation \eqref{eq: trace finite graph} in the case $t=1$ becomes
$$
M - \frac{1}{q}\mathrm{Tr}(A_X) = N_{0} + N_{1}q^{-1/2}\cdot \frac{-1}{\sqrt{q}}.
$$
Since $N_{0}=M$, we get that $N_{1} = \mathrm{Tr}(A_X)$.  In other words,  the number of closed irreducible walks of
length one equals the number of self-loops in the graph $X$, as expected.

The case when $t=2$ is the first instance where the second sum in Corollary \eqref{cor: trace fla}
is non-zero.  From \eqref{eq: I bessel as polynom} we
get that $I_{0}^{-2/\sqrt{q}}(2) = 1 + \frac{2}{q}.$
Additionally, \eqref{eq: I bessel as polynom} can be used to show that
$$
I_{1}^{-2/\sqrt{q}}(2) = \frac{-2}{2\sqrt{q}} = \frac{-2}{\sqrt{q}}
\,\,\,\,\,
\text{\rm and}
\,\,\,\,\,
I_{2}^{-2/\sqrt{q}}(2)
= \left(\frac{-2}{2\sqrt{q}}\right)^{2} = \frac{1}{q}.
$$
Therefore, from \eqref{eq: trace finite graph} we arrive at the expression
$$
M - \frac{2}{q}\mathrm{Tr}(A_X) + \frac{1}{q^{2}}\mathrm{Tr}(A^{2}_X) = N_{0}\left(1+\frac{2}{q}\right) + N_{1}\frac{-2}{q}
+ N_{2}\frac{1}{q^{2}} - \frac{M(q-1)}{q}.
$$
From the above evaluations of $N_{0}$ and $N_{1}$, we conclude that
\begin{equation}\label{eq: N2}
N_{2} = \mathrm{Tr}(A^{2}_X) - M(q+1).
\end{equation}

We can give a separate argument which proves \eqref{eq: N2} as follows.
It is well-known that $\mathrm{Tr}(A^{2}_X)$ denotes the number of closed walks of length $2$ on $X$; however,
$N_2$ is the number of closed walks of length $2$ that are irreducible and without tails. Closed walks of
length two with a tail are the same as closed reducible walks.  Hence, $N_2$ equals the difference between
$\mathrm{Tr}(A^{2}_X)$ and the number of closed, length $2$ reducible walks on $X$. For each vertex
$x \in V_X$ there are exactly $(q+1)$ reducible walks of length two; namely for each of the $q+1$
vertices $y_i$ adjacent to $x$, that is the walk $x,\, y_i,\, x$. Since the number of vertices is
$M$, this means that there are exactly $M(q+1)$ reducible walks of length two on $X$.  With all
this, one has a combinatorial argument which confirms \eqref{eq: N2}.

The calculations in the case $t=3$ follow a similar pattern as in the case $t=2$, but are a bit more involved.  From \eqref{eq: I bessel as polynom}
we obtain the following evaluations:
\begin{align*}
I_{0}^{-2/\sqrt{q}}(3) &= \sum\limits_{\alpha=0}^{1} \frac{3!}{\alpha!^2 (3-2\alpha)!}\left(\frac{-2}{2\sqrt{q}}\right)^{2\alpha}
 = 1 +\frac{6}{q};\\
I_{1}^{-2/\sqrt{q}}(3) &= \sum\limits_{\alpha=0}^{1} \frac{3!}{\alpha! (3-2\alpha-1)! (\alpha+1)!}\left(\frac{-2}{2\sqrt{q}}\right)^{2\alpha+1}
 = -\frac{3}{\sqrt{q}} -\frac{3}{q\sqrt{q}};\\
I_{2}^{-2/\sqrt{q}}(3) & = \left(\frac{-2}{2\sqrt{q}}\right)^{2} \cdot 3 = \frac{3}{q}; \quad \quad
I_{3}^{-2/\sqrt{q}}(3) = \left(\frac{-2}{2\sqrt{q}}\right)^{3} =\frac{-1}{q^{3/2}}.
\end{align*}
We now substitute these evaluations into \eqref{eq: trace finite graph} and get, upon expanding the left-hand-side, the
expression
\begin{align*}
M &- \frac{3}{q} \mathrm{Tr}(A_X) + \frac{3}{q^{2}} \mathrm{Tr}(A^{2}_X) - \frac{1}{q^{3}} \mathrm{Tr}(A^{3}_X)
\\&= N_{0}\left(1+\frac{6}{q}\right) -\frac{1}{q} N_{1}\left(3+\frac{3}{q}\right)+ N_{2}q^{-1}\cdot \frac{3}{q}
+ N_{3}q^{-3/2}\cdot \frac{-1}{q^{3/2}} - M(q-1)q^{-1}\cdot \frac{3}{q}.
\end{align*}
When solving for $N_{3}$, one uses the values for $N_{0}$, $N_{1}$ and $N_{2}$ obtained above to get that
\begin{equation}\label{eq: N3}
N_{3} = \mathrm{Tr}(A^{3}_X) -3q \mathrm{Tr}(A_X).
\end{equation}

Let us provide a combinatorial proof of formula \eqref{eq: N3}.
By definition, $N_3$ is the number of closed, irreducible walks of length $3$
on $X$ without tails. As such, $N_{3}$ equals the difference between $\mathrm{Tr}(A^{3}_X)$,
the number of all closed walks of length $3$ on $X$ and the number of closed
walks of length $3$ that are reducible or possess a tail. For vertices $x\in V_X$
which do not possess self-loops, all closed walks of length $3$ are irreducible
and without tails. Hence, if there are no self-loops in $X$, then $\mathrm{Tr}(A_X)=0$
and the formula \eqref{eq: N3} holds true. Assume now that there are $N_1>0$ vertices
in $V_X$ which are adjacent to themselves, meaning they have self-loops. Let $x$ be such a vertex and let $y\neq x$
be a vertex adjacent to $x$; there are $q$ such vertices. Then, the walks $x,\,x,\,y,\,x$
and $x,\,y,\,x,\,x$ are all reducible walks starting from $x$ and for each such $x$
there are $2q$ of them. Moreover, the walks $y,\,x,\,x,\,y$ are only possible length
three walks that possess a tail; for a fixed $x$ with a self-loop, there are exactly
$q$ of them. This proves that the number of closed walks on $X$ of length $3$ which
are either reducible or possess a tail equals $3qN_1=3q \mathrm{Tr}(A_X)$ and completes
the proof of \eqref{eq: N3}.

The formulas for $N_{j}$ for $j \leq 3$ obtained in these calculations agree with the expressions obtained from equation (34) of \cite{Mn07}.


\subsection{Evaluating $N_{j}$ for general $j$}\label{sec: evaluating_N_m}
Going further, one can use \eqref{eq: I bessel as polynom} to obtain closed-form expressions
for each $\tilde{N}_t:=q^{-t/2}N_{t}$ for any finite $(q+1)$-regular graph with $M$ vertices.  The
approach is as follows.

Fix $t\geq 1$ and let $V_{t}$ denote the $(t+1)\times (t+1)$ matrix whose $(j,k)$
entry is $v_{j,k}=I_{k-1}^{-2/\sqrt{q}}(j-1)$ for $j,k=1,\ldots,t+1$.  Observe that $v_{j,k}= 0$ for $k > j$; see Lemma 4(i).  Hence,
$V$ is a lower triangular matrix.  Furthermore, from \eqref{eq: I bessel as polynom} we get that
$$
v_{k,k} = I_{k-1}^{-2/\sqrt{q}}(k-1) = (-1)^{k-1} q^{-(k-1)/2}
\,\,\,\,\,
\text{\rm for}
\,\,\,\,\,
k=1,\ldots, t+1.
$$
Let $\mathbb{T}$ be the column vector of length $t+1$ whose $k$-th entry is
$\mathrm{Tr}\left(\left(\frac{A_X}{2\sqrt{q}}\right)^{k-1}\right)$, where the first entry is $\mathrm{Tr}(\mathrm{Id})=M$.
Let $\tilde{\mathbb{N}}$ be the column vector of length $t+1$ whose $k$-th entry is $\tilde{N}_{k-1}$.
Finally, let $\mathbb{E}=(e_k)_{k=1}^{t+1}$ be the column vector of length $t+1$ whose $k$-th entry is
$M(q-1)q^{-1}I^{-2/\sqrt{q}}_{2}(k-1)$.
With all this, the set of equations \eqref{eq: trace finite graph} can be written as
$$
B\mathbb{T} = V \tilde{\mathbb{N}} - \mathbb{E},
$$
where $B=(b_{jk})_{(t+1)\times (t+1)}$ is the matrix with entries $b_{jk}=\binom{k-1}{j-1}
\left(-\frac{2}{\sqrt{q}}\right)^{j-1}$, for $k\geq j$ and $b_{jk}=0$ for $k<j$.  By convention we set $b_{11}=\binom{0}{0}=1$.

The matrix $V$ is invertible, hence the solution to the set of equations  \eqref{eq: trace finite graph} for any fixed $t\geq 1$ is
\begin{equation}\label{eq: tilde matrix}
\tilde{\mathbb{N}}=V^{-1}\left(B\mathbb{T} + \mathbb{E}\right).
\end{equation}

Let $V^{-1}=(\tilde{v}_{jk})_{(t+1)\times (t+1)}$.  Let us show how to compute $V^{-1}$
in terms of powers of discrete $I$-Bessel functions. Let $D$ be the $(t+1)\times (t+1)$ diagonal matrix whose $(k,k)$ entry is
 $d_{k,k} =(-1)^{k-1} (q)^{-(k-1)/2}$.  Then $D^{-1}$ is the diagonal matrix whose
$k$-th diagonal element is $(-1)^{k-1} (q)^{(k-1)/2}$.  Let us write $V = D + \tilde{V}$, so $\tilde{V}$ is a strictly lower diagonal
matrix.

Since $\tilde{V}$ is strictly lower triangular and $D^{-1}$ is diagonal, then $D^{-1}\tilde{V}$ is also strictly lower triangular,
so then $(D^{-1}\tilde{V})^{t+1}$ is the zero matrix.  Therefore,
\begin{equation}\label{eq:V_inverse}
V^{-1} = (I+D^{-1}\tilde{V})^{-1}D^{-1}=
\left(\sum\limits_{h=0}^{t}\left(-D^{-1}\tilde{V}\right)^{h}\right)D^{-1},
\end{equation}
which shows that entries $\tilde{v}_{jk}$ of $V^{-1}$ can be written in terms of sums of products of discrete $I$-Bessel functions. Therefore, formula \eqref{eq: tilde matrix} combined with \eqref{eq:V_inverse} provides an efficient algorithm for evaluation of $\tilde{N}_j$ for general $j$ in terms of values of discrete $I-$Bessel functions.

\subsection{Chebyshev polynomials and discrete $I$-Bessel functions}

In this section we will compare our results with equation (34) of \cite{Mn07}.
Recall that \cite{Mn07} assumed $X$ is a finite connected $(q+1)$-regular graph with $M$ vertices, without multiple edges and loops.
So, for this section we also will assume that $X$ does not have multiple edges or self-loops.

For any positive integer $\ell$, Equation (34) of \cite{Mn07} can be written, in our notation, as
\begin{equation}\label{eq:Mnev_counting}
  \tilde{N}_{\ell}=2 \mathrm{Tr}\left(T_{\ell}\left(\frac{A_X}{2\sqrt{q}}\right)\right)+ \frac{1+(-1)^{\ell}}{2}(q-1)q^{-\ell/2}M,
\end{equation}
where $T_{\ell}(x)$ denotes the $\ell$-th Chebyshev polynomial of the first kind. In other
words, results from \cite{Mn07} evaluate the left-hand side of \eqref{eq: tilde matrix} in
terms of the Chebyshev polynomial; this evaluation can be viewed as a linear form in variables
$$
x_j(X,q):=\mathrm{Tr}\left(\left(\frac{A_X}{2\sqrt{q}}\right)^j\right)
\,\,\,\,\,
\text{\rm for $j=0,1,\ldots \ell$,}
$$
together with the term $\frac{1+(-1)^{\ell}}{2}(q-1)q^{-\ell/2}M$,
which is zero if $\ell$ is odd. Recall that we adopt the convention that
$$
x_0(X,q):=\mathrm{Tr}
\left(\left(\frac{A_X}{2\sqrt{q}}\right)^0\right)=M.
$$

In other words, we can write equation \eqref{eq:Mnev_counting} as
\begin{equation}\label{eq:Mnev_counting 1}
  \tilde{N}_{\ell}=2 \sum_{j=0}^{\ell} t_{\ell}(j) x_j(X,q)  + \frac{1+(-1)^{\ell}}{2}(q-1)q^{-\ell/2}M,
\end{equation}
where $t_{\ell}(j)$ is the coefficient multiplying $y^j$ in the expansion of the Chebyshev polynomial $T_{\ell}(y)$.

Returning to our calculations, note that the $k$-th entry on the right-hand side of \eqref{eq: tilde matrix}
equals a linear form in the variable $x_j(X,q)$, $j=0,\ldots,k-1$, namely
$$
P_{k-1}(x_X(q))=\sum_{j=0}^{k-1} a_{k-1,j}x_j(X,q),
$$
where $a_{k-1,j}=\sum_{m=1}^{k} \tilde{v}_{km}b_{m(j+1)}$, for $j=1,\ldots, k-1$ and $a_{k-1,0}=\sum_{m=1}^{k} \tilde{v}_{km}(b_{m1} + e_{k-1})$.

Therefore, equation \eqref{eq: tilde matrix} combines with \eqref{eq:Mnev_counting 1} to imply
that for any $k\geq 0$ we have the expression
$$
2 \sum_{j=0}^{k} t_{k}(j) x_j(X,q)+\frac{1+(-1)^{k}}{2}(q-1)q^{-k/2}M = \sum_{j=0}^{k} a_{k,j}x_j(X,q).
$$
The coefficients of the linear form on the right-hand side are explicitly computable in terms of discrete $I$-Bessel functions
$I_{j}^{-2/\sqrt{q}}(\ell)$ for $j\leq \ell \leq k-1$. In other words, the rows of the matrix
$$
\left(\sum\limits_{h=0}^{t}\left(-D^{-1}\tilde{V}\right)^{h}\right)D^{-1}B
$$
amount a listing of the coefficients of the Chebyshev polynomial associated to non-constant terms.

We find it intriguing to see that the special values of the discrete $I$-Bessel function come together
through the above computations and produce the coefficients of Chebyshev polynomials.

\subsection{Further counting problems}\label{sec: further counting}

The analysis of the previous sections provided an explicit, closed-form expression for $N_{m}$, the number
of closed irreducible walks without tails of length $m$ on a finite $(q+1)$-regular graph, in terms of a finite sum involving
the length spectrum of the adjacency matrix.  In the setting of continuous time heat kernels, as in \cite{CJK14}
or \cite{Mn07}, the authors obtain expressions which relate a spectral term, often finite, to a geometric
term, which is an infinite sum in the continuous time setting.  In the discrete time setting we study in
this article, the geometric side of the identities we study also is finite.  As a result, the formula for
$N_{m}$ follows from our heat kernel formulas and elementary linear algebra.

Let us now briefly describe how one can expand the approach taken above to the case when $X$ is infinite.

One begins with the pre-trace formula as given in Corollary \eqref{cor: pre-trace}.  It is not necessary
to take $x=x_{0}$, but certainly that is a possibility.  For a fixed integer $t$, the left-hand side of
the formulas in Corollary \eqref{cor: pre-trace} could be taken as written or expanded into a series involving
the various moments of the spectral measure.  As in Section \ref{sec: evaluating_N_m}, one can solve for
$b_{m}(x)$, ultimately obtaining an expression for $b_{t}(x)$ in terms of the $m$-th moments of the spectral
measure, for all $m \leq t$, and special values of the discrete $I$-Bessel function.  Equivalently, the formulas
for $b_{t}(x)$ will involve the shifted moments
$$
\int\limits_{-(q+1)}^{(q+1)}(\lambda - q)^{m} \mu_{x}(d\lambda)
\,\,\,\,\,
\text{\rm for $m \leq t$}
$$
and special values of the discrete $I$-Bessel function.  The resulting formulas
will be similar to those in \eqref{eq: tilde matrix} and \eqref{eq:V_inverse}.
Finally, one can go one step further and get expressions for $c_{m}(x)$ using
Lemma \ref{lem:data_relations} and one further matrix multiplication.

\section{Applications to other diffusion processes} \label{sec: applications}

Let us briefly indicate how to apply our main results to give a closed formula for the solution of a certain zero-sum diffusion process on the $(q+1)$-regular tree $T_{q+1}$, which can be viewed as a random walk in presence of a heat or cooling source.

Namely, let $\alpha,\beta >0$ be such that $\alpha + \beta<1$.  Consider the zero-sum diffusion process on the non-negative integers
$\mathbb{N}_0$ as described by the difference equation
\begin{equation}\label{eq:random walk on half line}
 \partial_t K(x;t)= \left\{
                       \begin{array}{ll}
                         (\beta-1) K(0;t) + (1-\beta)K(1;t), & \text{  if  } x=0; \\
                         (\beta-1) K(x;t)+(1-\beta-\alpha)K(x+1;t)+ \alpha K(x-1;t), & \text{  if  } x\geq 1,
                       \end{array}
                     \right.
\end{equation}
with initial condition as in \eqref{eq: initial cond tree}.  One can view the conditions at $x=0$ as that of a reflecting
boundary.

Equation \eqref{eq: heat eq on a tree} can be suitably rescaled to give a solution to \eqref{eq:random walk on half line}.
Indeed, by reasoning as in the proof of Proposition \ref{Prop: HK on tree},  one can show that
the solution $K(x;t)$ for $x,t\in \mathbb{N}_0$  is given by
$$
K(x;t)=\beta^t\left[\left(\frac{1-\alpha-\beta}{\alpha}\right)^{-x/2}I_x^{a}(t) + \frac{2\alpha+\beta-1}{\alpha}\sum_{j=1}^{\infty} \left(\frac{1-\alpha-\beta}{\alpha}\right)^{-(x+2j)/2}I_{x+2j}^{a}(t)\right],
$$
where $a=2\sqrt{\frac{\alpha(1-\alpha-\beta)}{\beta^2}}$.  We will leave the proof to the interested reader.

\vspace{5mm}\noindent
Carlos A. Cadavid \\
Department of Mathematics \\
Universidad Eafit \\
Carrera 49 No 7 Sur-50 \\
Medell\'in, Colombia \\
e-mail: ccadavid@eafit.edu.co

\vspace{5mm}\noindent
Paulina Hoyos  \\
Department of Mathematics \\
The University of Texas at Austin \\
C2515 Speedway, PMA 8.100 \\
Austin, TX 78712
U.S.A. \\
e-mail: paulinah@utexas.edu

\vspace{5mm}\noindent
Jay Jorgenson \\
Department of Mathematics \\
The City College of New York \\
Convent Avenue at 138th Street \\
New York, NY 10031
U.S.A. \\
e-mail: jjorgenson@mindspring.com

\vspace{5mm}\noindent
Lejla Smajlovi\'c \\
Department of Mathematics \\
University of Sarajevo\\
Zmaja od Bosne 35, 71 000 Sarajevo\\
Bosnia and Herzegovina\\
e-mail: lejlas@pmf.unsa.ba

\vspace{5mm}\noindent
Juan D. V\'elez \\
Department of Mathematics \\
Universidad Nacional de Colombia\\
Carrera 65 Nro. 59A - 110\\
Medell\'in, Colombia\\
e-mail: jdvelez@unal.edu.co

\end{document}